\documentclass{m2an}
\usepackage[utf8]{inputenc}
\usepackage{mathtools,cases,csquotes}
\usepackage{amsmath}
\usepackage{cancel}
\usepackage{amsbsy,csquotes}
\usepackage{booktabs}
\usepackage{subcaption}
\usepackage{amsmath, amssymb}
 \usepackage{listings}
\usepackage{empheq}
\usepackage{afterpage}

\DeclareCaptionLabelFormat{andtable}{\tablename \hspace{0.1em} #2}

\newcommand{\norm}[2]{ \Vert #1 \Vert_{#2}}

\begin{document}
\title{Splitting method for elliptic equations with line sources}
\thanks{This work was partially supported by the Research Council of Norway,  project number 250223, and Deutsche Forschungsgemeinschaft, grant number WO-671/11-1.}
\author{Ingeborg G. Gjerde}\address{Department of Mathematics, University of Bergen, Norway. \email{ingeborg.gjerde@uib.no, jan.nordbotten@uib.no}}
\author{Kundan Kumar}\address{Department of Mathematics, Karlstad University, Sweden.  \email{kundan.kumar@kau.se}}
\author{Jan M. Nordbotten}\sameaddress{1}\secondaddress{Department of Civil and Environmental Engineering, Princeton University, USA.}
\author{Barbara Wohlmuth}\address{Department of Mathematics, Technical University of Munich, Germany. \email{ barbara.wohlmuth@ma.tum.de}}
\date{October 25, 2018}
\begin{abstract}
In this paper, we study the mathematical structure and numerical approximation of elliptic problems posed in a (3D) domain $\Omega$ when the right-hand side is a (1D) line source $\Lambda$. The analysis and approximation of such problems is known to be non-standard as the line source causes the solution to be singular. Our main result is a splitting theorem for the solution; we show that the solution admits a split into an explicit, low regularity term capturing the singularity, and a high-regularity correction term $w$ being the solution of a suitable elliptic equation. The splitting theorem states the mathematical structure of the solution; in particular, we find that the solution has anisotropic regularity. More precisely, the solution fails to belong to $H^1$ in the neighbourhood of $\Lambda$, but exhibits piecewise $H^2$-regularity parallel to $\Lambda$. The splitting theorem can further be used to formulate a numerical method in which the solution is approximated via its correction function $w$. This approach has several benefits. Firstly, it recasts the problem as a 3D elliptic problem with a 3D right-hand side belonging to $L^2$, a problem for which the discretizations and solvers are readily available. Secondly, it makes the numerical approximation independent of the discretization of $\Lambda$; thirdly, it improves the approximation properties of the numerical method. We consider here the Galerkin finite element method, and show that the singularity subtraction then recovers optimal convergence rates on uniform meshes, i.e., without needing to refine the mesh around each line segment. The numerical method presented in this paper is therefore well-suited for applications involving a large number of line segments. We illustrate this by treating a dataset (consisting of $\sim 3000$ line segments) describing the vascular system of the brain.
\end{abstract}

 \maketitle

\section{Introduction}

\def\avint{\mathop{\,\rlap{-}\!\!\int}\nolimits}

Let $\Omega \subset \mathbb{R}^3$ be an open domain with the smooth boundary $\partial \Omega$, $\Lambda = \cup_{i=1}^n \Lambda_i$ be a collection of line segments $\Lambda_i$, and $\kappa \in W^{2, \infty}(\Omega)$ be a uniformly positive, scalar-valued coefficient. We consider elliptic problems with line sources of the type
\begin{subequations}
\begin{align}
 - \nabla \cdot (\kappa \nabla u) &= f \delta_\Lambda \quad   \text{in} \, \Omega, \label{eq:basic1}\\
 \int_\Omega f \delta_\Lambda \, v \, \mathrm{d} \Omega  &= \sum_{i=1}^n \int_{\Lambda_i} f(s_i) v(s_i) \mathrm{d}s \quad \forall v \in C^0(\Omega),
 \label{eq:basic2}
\end{align}
\end{subequations}
$s_i$ being the arc-length of line segment $i$. The right-hand side $f$ is supported on $\Lambda$ and models a collection of 1D sources and sinks. The line source can be interpreted as a 1D fracture in the domain, with the variable $f : \Lambda \rightarrow \mathbb{R}$ denoting the linear mass flux from $\Lambda$ into the domain. This fracture is modelled mathematically by means of a Dirac distribution on the line, understood in the sense of \eqref{eq:basic2}.

Models of this type arise in a variety of applications, e.g., in the modelling of 1D steel components in concrete structures \cite{llau2016} or the interference of metallic pipelines and bore-casings in electromagnetic modelling of reservoirs \cite{weiss2017}. Of special interest are coupled 1D-3D models, where \eqref{eq:basic1} is coupled with a 1D flow equation defined on $\Lambda$ \cite{dangelo2008}. The coupled 1D-3D model is more commonly used for biological applications, e.g., in the study of blood flow in the  vascularized tissue of the brain \cite{Reichold2009}, \cite{Grinberg2011}, the efficiency of cancer treatment by hyperthermia \cite{nabil2016} and by drug delivery through microcirculation \cite{Cattaneo2014}, \cite{zunino2018}. The use of a 1D source term is motivated by the observation that roots, blood vessels, steel components and pipelines all have negligible radii compared to the size of the simulation domain. It would therefore be quite expensive to resolve these structures as 3D objects in a mesh, and we model them instead as 1D structures embedded in $\Omega$. 

The coupled 1D-3D problem is a particular instance of mixed-dimensional PDEs, where an $d$-dimensional PDE is coupled to some $d-k, k < d$-dimensional PDE. A typical example can be found in fracture modelling, where a fracture is considered as a $d-1$ dimensional manifold embedded in a $d$-dimensional domain. We refer to \cite{Fracture2005SISC}, \cite{ Alessiofracture}, \cite{Kumarfracture} and references therein for dealing with some examples, \cite{miro2016-2D1D} for suitable preconditioners, and \cite{boon2017} for a general mathematical framework  to handle the $d$ to $d-1$-dimensional coupling. The results in these works are, however, not easily transferred to the coupled 1D-3D problem, as the dimensional gap impacts the mathematical structure of the solution. In particular, an increase in the dimensional gap negatively affects the solution regularity. In fact, the line source in the right-hand side of \eqref{eq:basic1} induces a logarithmic type singularity in the solution, ultimately causing the solution not to be in $H^1(\Omega)$. The approximation and analysis of the line source problem therefore requires special consideration. 

\begin{figure}
\centering
\includegraphics[width=0.3\textwidth]{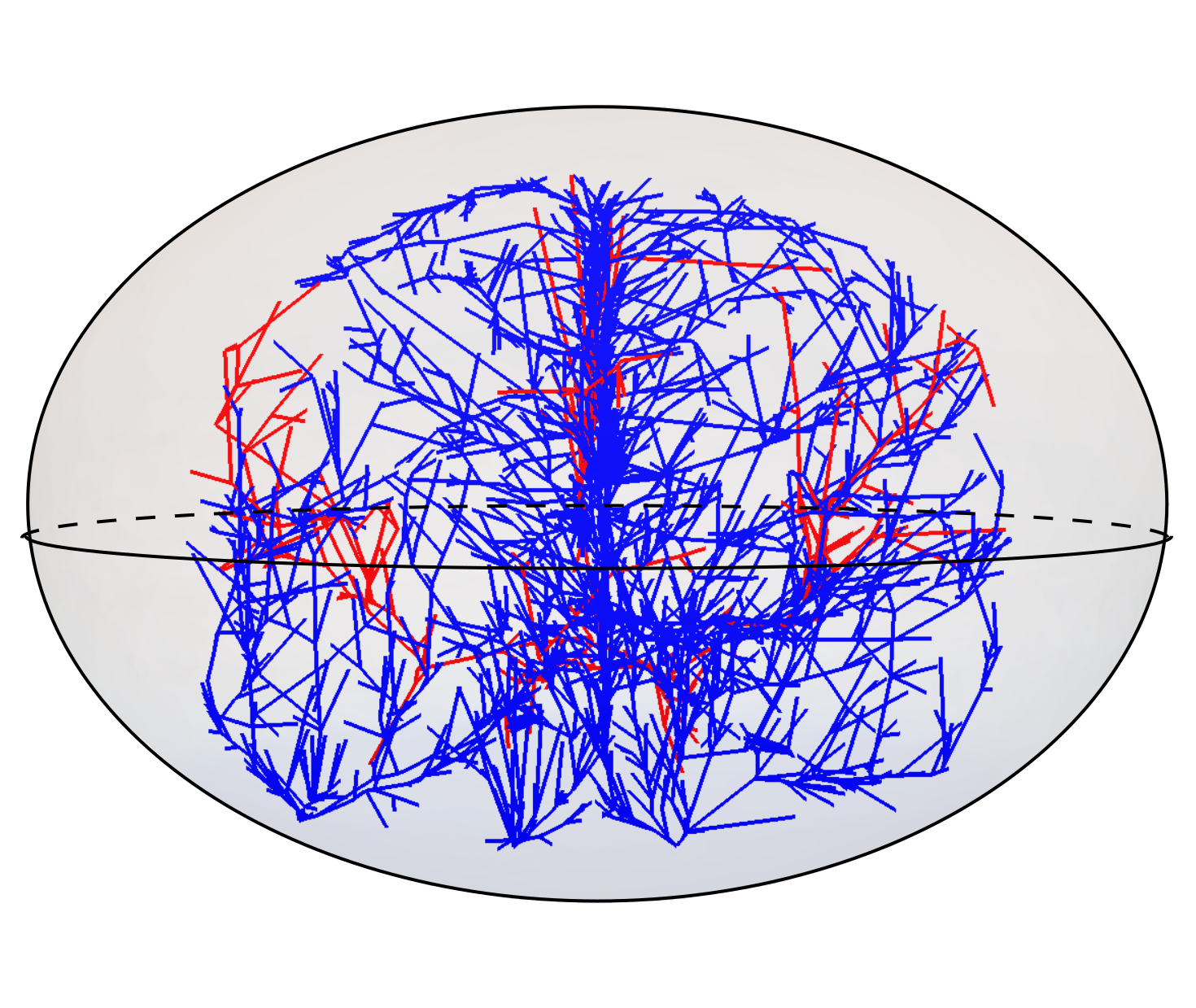}
\includegraphics[width=0.3\textwidth]{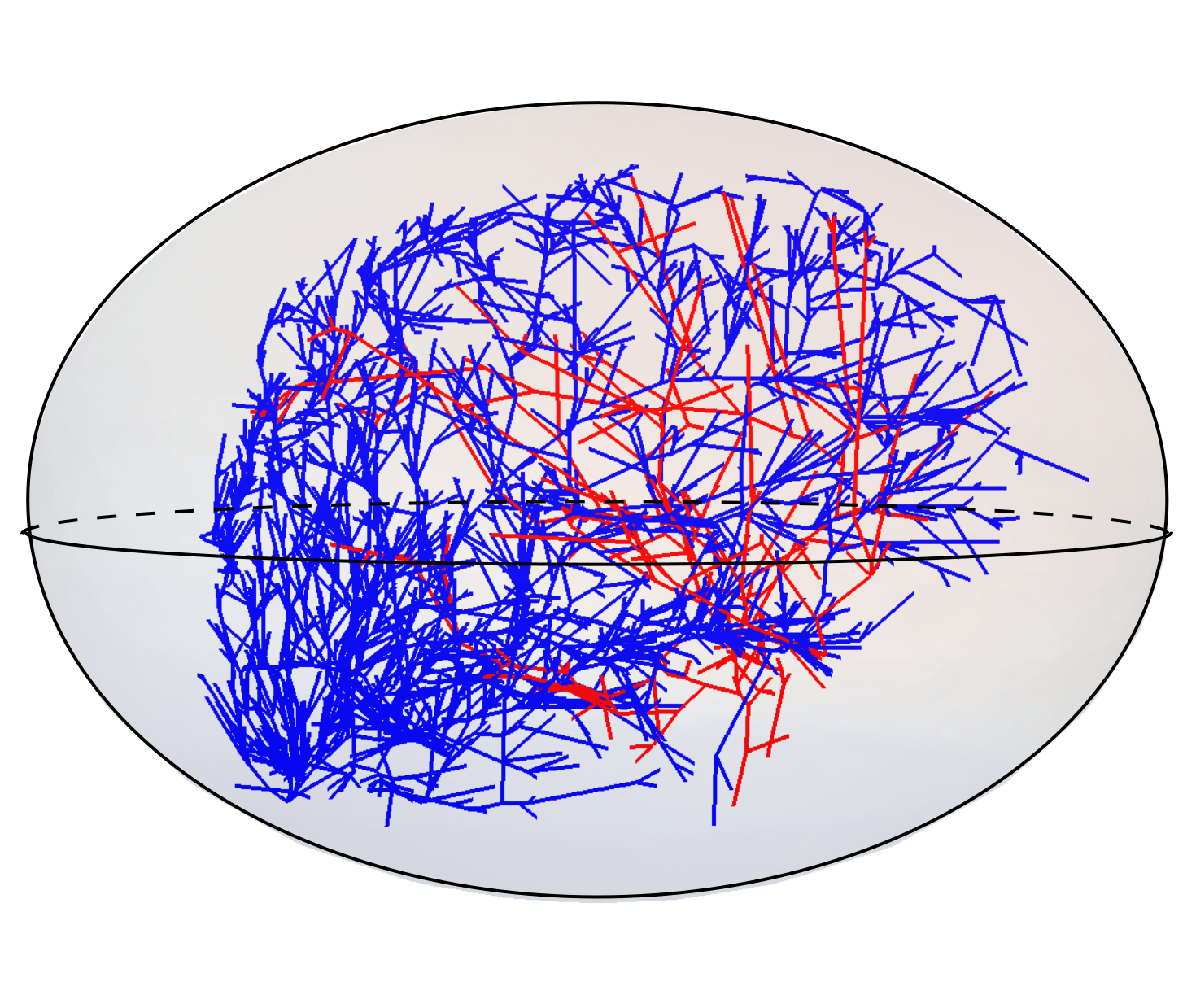}
\includegraphics[width=0.3\textwidth]{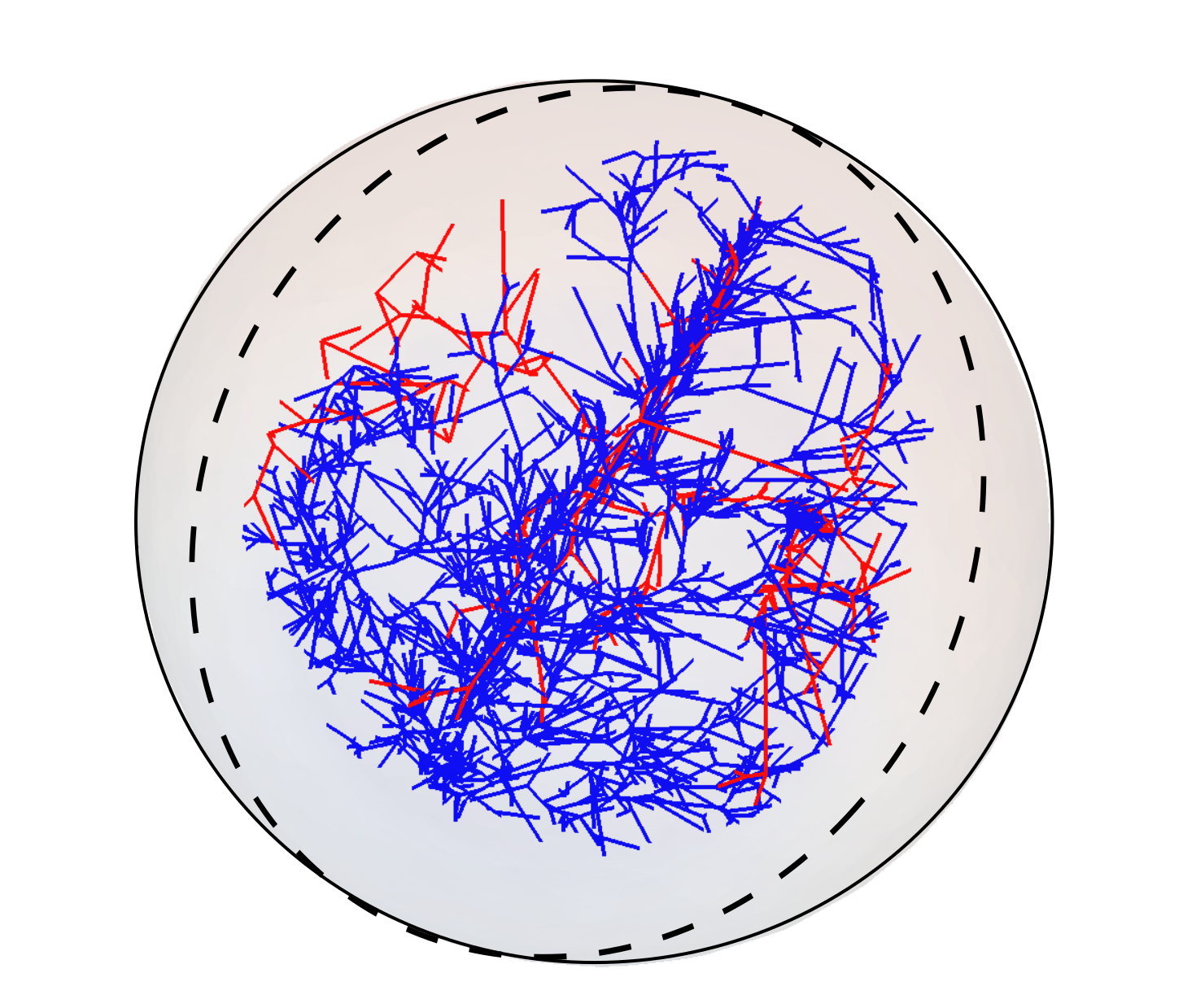}
\caption{Front, side and perspective view of the arteries (red) and veins (blue) in a human brain. The arterial system is responsible for providing oxygenated blood to the brain tissue; its counterpart, the venous system, then carries deoxygenated blood back to the heart. Both systems are composed of a collection of blood vessels that individually resemble line segments $\Lambda_i$. It is therefore natural to model the arterial and venous systems as graphs $\Lambda_a =\cup_{i=1}^n \Lambda_{a,i}$ and $\Lambda_v=\cup_{i=1}^n \Lambda_{v,i}$ respectively. The dataset consists of $\sim 3000$ line segments and was constructed from an MRI scan \cite{brain}.}
\label{fig:drawing}
\end{figure}

The existence of a solution to \eqref{eq:basic1} was proven in \cite{dangelo2008} by means of weighted Sobolev spaces. Convergence rates for the Galerkin finite element (FE) method were proven in \cite{dangelo2012} in terms of weighted Sobolev spaces, with results similar to those found in \cite{apel2011} for the 2D problem with a point source. In both of these works, the high dimensional gap was found to cause sub-optimal convergence rates on uniform meshes, with the $L^2$-error converging with order less than one, and no convergence in the $H^1$ norm. These sub-optimal convergence rates for the 2D-0D problem were found to be local to the point source in \cite{koppl2017}, \cite{koppl2015}, in the sense that linear finite elements converge \ optimal up to a log-factor in the $L^2$-norm if a small neighbourhood is removed around the point. A similar result was proved for the 3D problem with a line source in \cite{koppl2017}. Optimal convergence rates can be retrieved in the entire domain by use of a graded mesh, i.e., by a particular refinement of the mesh around each point \cite{apel2011} or line segment \cite{dangelo2012}. 

The numerical approximation of coupled 1D-3D models is similarly known to be difficult due to the presence of a line source in the 3D flow equation. We refer to \cite{miro2016-2D1D} for preconditioners developed to treat the ill-conditioning this coupling introduces in the discretization matrix $A$. In this work, negative fractional Sobolev norms were found suitable to precondition the Schur complement of the system, with an extension to the positive fractional Sobolev norms considered in \cite{baerland2018}. In \cite{koppl2017}, an alternative coupling scheme was introduced, where the source term was taken to live on the boundary of the inclusions. The result is a 1D-(2D)-3D method where the dimensional gap has been reduced to 1. This improves the regularity and thus the approximation properties of the problem, at the expense of having to resolve the 2D boundary of the inclusions.

Our approach is based on exploiting certain splitting properties of the solution, where the solution $u$ is split into a collection of low-regularity logarithmic terms $G_i$, each corresponding to a line segment $\Lambda_i$, and a high-regularity correction term $w \in H^2(\Omega)$:
\begin{align}
u = \sum_{i=1}^n E_i(f) G_i + w.
\label{eq:decomp1}
\end{align}
The term $E_i(f)$ here denotes an extension of $f$ from the line segment $\Lambda_i$ to the entire domain, the precise definition of which will be given at a later point. The logarithmic term was found by an integration of the 3D Green's function for the Laplacian over each line segment, resulting in an analytic function for the line source similar to the one derived in \cite{strack2017analytical}. The correction term $w$ solves its own elliptic equation on $\Omega$ with a given right-hand side $F \in L^{2-\epsilon}(\Omega)$ for some arbitrarily small $\epsilon>0$. $F$ depends on the line source intensities $f_i$, the endpoints of each line segment, and the coefficient $\kappa$; its regularity can be improved to $F \in L^2(\Omega)$ given some assumptions on the structure of $\kappa$.

The $G_i$ terms completely capture the singular parts of the solution, and the sum $\sum_{i=1}^n f_i G_i$ therefore constitutes the low-regularity part of the split. The mathematical structure and regularity properties of $u$ can thus be found by inspection of this term. In particular, this will reveal to us that the solution $u$ exhibits anisotropic regularity, with the low regularity of $u$ occurring in the plane perpendicular to the line. Formally, the line singularity acts at each surface normal to a line segment as a Dirac source term in a 2D plane, whereas along the line it inherits the regularity of $f$. 

Besides informing us about the mathematical structure of the solution, the splitting in \eqref{eq:decomp1} suggests an alternative numerical method in which $u$ is solved for indirectly by a correction term $w$. A similar approach is known for the point source problem, where it is often referred to as a singularity subtraction technique \cite{Drechsler2009}. With this in mind, we will refer to our method as the \textit{singularity subtraction technique for the line source problem.} The singularity subtraction recasts \eqref{eq:basic1}-\eqref{eq:basic2} as a (3D) elliptic problem with a (3D) right-hand side belonging to $L^2(\Omega)$, a standard problem for which suitable approximation methods are well known. As the splitting works at the continuous level, this method is readily adapted to different discretization schemes; for example, we provide in \cite{ecmor} a mixed method adaptation of the ideas developed here.  Furthermore, the technique can be adapted to handle also coupled 1D-3D problems. By a decoupling, we can solve the coupled 1D-3D problem iteratively, solving at each iteration the 3D elliptic problem \eqref{eq:basic1} with a \textit{given }line source intensity $f$.  We provide in \cite{ecmor} numerical evidence that this iterative solution method converges when the singularity subtraction technique is applied. For a convergence proof of a similar iterative scheme, we refer here to \cite{kopplphd}.

The singularity subtraction is associated with a number of advantages, the main ones being that it results in numerical methods that are (i) easy to implement, (ii) can be adapted to a variety of discretization methods, and (iii) have significantly improved convergence properties. We consider here the use of a Galerkin FE method, in which case the singularity subtraction technique can be considered as a special case of the Generalized Finite Element Method \cite{gfem}. The Generalized Finite Element Method works by enriching the basis functions at locations where the analytical solution corresponding to the singularity is known (for example, displacement in a corner in an elasticity equation). This improves the convergence properties of the method; we show in this work that the singularity subtraction recovers optimal convergence rates on uniform meshes. 

As the singularity subtraction makes the problem independent of the discretization of the lines, the numerical method presented in this paper is particularly well suited for problems involving a large number of line sources. This type of problem arises for example when considering flow through the vascularized tissue of the brain; the dataset describing the arterial and venous systems of the brain (shown in Figure \ref{fig:drawing}) consists of $\sim 3000$ line segments. Moreover, this dataset is incomplete in the sense that the imaging captures only blood vessels of a certain radii. The size of this type of dataset can therefore be expected to increase as image quality is increased with improved technology. The main limitation of our approach is that it handles only a certain collection of problems. The study presented here requires the problem to be linear, moreover, the permeability coefficients has to be scalar-valued and smooth. This can be relaxed to have locally smooth coefficients, however, when the number of line segments is large, ensuring a smooth permeability even locally is restrictive. 

Our approach is similar to the seminal 
work of Peaceman \cite{peaceman1978well} in that way it uses the properties of the analytical solution to inform and improve the approximation of more general problems. The Peaceman well correction is used to account for the fact that bottom-hole pressure in a well differs greatly from the cell pressure of the cell in which it is located \cite{Ewing1999}. The correction factor itself is strongly dependent on the spatial relationsship between well and discretization, thus requiring various adaptations in order to handle different types of grids and discretizations \cite{Lie2012}, \cite{aavatsmark2003well}. The singularity subtraction, conversely, works on the continuous level, thus making it independent of the specific discretization method.

The article is structured as follows. In Section \ref{sec:notation}, we introduce the notation that will be used to state the regularity properties of the solution. We introduce here both weighted and fractional Sobolev spaces, as these are the spaces we have seen used for the coupled 1D-3D problem in the aforementioned references. In Section \ref{sec:decomposition}, we show that the solution admits a splitting into high and low-regularity terms, and discuss how this informs us of the mathematical structure of the solution. We start in Section \ref{sec:infinite-line} with the special case of $f$ smooth and $\Lambda$ consisting of a single line segment. With these assumptions the splitting is quite trivial, making it easier to illustrate the splitting technique itself and discuss the regularity properties of the solution. We then proceed to generalize the splitting technique, showing in Section \ref{sec:finite-line} how to treat an arbitrary line segment $\Lambda_i$, and then in Section \ref{sec:kappa} how to treat the case of $\kappa$ varying spatially. The decomposition section is then concluded in Section \ref{sec:existence} with a splitting theorem valid for any graph $\Lambda$.

The next two sections are dedicated to the numerical approximation of \eqref{eq:basic1}; from this point on, we treat only the case of $\kappa=1$. In Section \ref{sec:numerical-disc}, we present two different FE methods, the first being the straight-forward discretization of \eqref{eq:basic1}-\eqref{eq:basic2}. In this method, the line source is assembled explicitly into the right-hand side as a line integral; this is the same discretization as was presented in \cite{dangelo2012}. The second method is a singularity subtraction approach based on the splitting theorem in Section \ref{sec:existence}, where the correction function $w$ is approximated using FE. We denote this method as the singularity subtraction based finite element (SSB-FE) method. In Section \ref{sec:numerical-disc}, we then show how the singularity subtraction leads to an improvement in the convergence rates of the method. We also discuss the modelling error introduced by neglecting some of the line segments in $\Lambda$. The article is then concluded with three numerical experiments. The first two of these serve to verify the mathematical properties and error estimates found in Sections \ref{sec:decomposition} and \ref{sec:numerical-disc}, respectively. The last example uses the graph shown in Figure \ref{fig:drawing} to demonstrate the capabilities of the SSB-FE method in handling datasets with a large number of line segments. 

\section{Background and notation}
\label{sec:notation}
In this section we will introduce the function spaces used in this work. Let $W^{k, p}(\Omega)$ be the standard Sobolev space, 
\begin{align*}
W^{k, p}(\Omega) = \lbrace u \in L^p(\Omega) : D^\beta u \in L^p(\Omega) \text{ for } \vert \beta \vert \leq k \rbrace,
\end{align*}
with $\beta$ denoting a multi-index and $D^\beta$ the corresponding weak distributional derivative of $u$. This is a Banach space endowed with the norm
\begin{align*}
\norm{u}{W^{k, p}(\Omega)}^p = \sum_{\vert \beta \vert < k} \norm{D^\beta u}{L^p(\Omega)}^p.
\end{align*}
For $p=2$ it is also a Hilbert space, customarily denoted $H^k(\Omega)=W^{k,2}(\Omega)$, with inner product
\begin{align*}
(u, v)_{H^k(\Omega)} = \sum_{\vert \beta \vert < k} \int_\Omega D^{\beta} u \, D^{\beta} v\, \mathrm{d} \Omega. 
\end{align*}
As we will see, solutions to \eqref{eq:basic1}-\eqref{eq:basic2} fail to belong to $H^1(\Omega)$; for this reason, we consider also weighted and fractional Sobolev spaces. To define the weighted Sobolev space, let $ -1 < \alpha < 1$, and take $L^2_\alpha$ to denote the weighted Hilbert space
\begin{align*}
L_\alpha^2(\Omega) := \left \{u: \int_\Omega u^2 r^{2\alpha} \mathrm{d} \Omega < \infty \right \},
\end{align*}
equipped with the inner product 
\begin{align*}
(u, v)_{L_\alpha^2} = \int_\Omega r^{2\alpha} u v \, \mathrm{d} \Omega.
\end{align*}
For $\alpha>0$, the weight $r^\alpha$ has the power to dampen out singular behaviour in the function being normed, and we have the relation $L^2_{-\alpha}(\Omega) \subset L^2(\Omega) \subset L^2_\alpha(\Omega)$. A weighted Sobolev space $V^k_\alpha(\Omega)$ can now be constructed as \cite{koslov}, \cite{oleinik}
\begin{align*}
V^k_\alpha(\Omega) = \lbrace D^\beta u \in L_{\alpha + \vert \beta \vert - k}^2(\Omega) , \vert \beta \vert \leq m \rbrace,
\end{align*}
equipped with the norm
\begin{align*}
 \norm{u}{V^k_\alpha(\Omega)}^2 = \sum_{\vert \beta \vert \leq k} \Vert D^\beta u \Vert_{L^2_{\alpha+\vert \beta \vert - k}}^2.
\end{align*}
This weighted space takes the order of the derivative into consideration when computing the weight function, increasing the weighting-factor to compensate for the regularity loss due to an extra derivative. We shall use mainly ${V^1_\alpha(\Omega)}$, in which the norm can be stated more simply as $\norm{u}{V^1_\alpha(\Omega)}^2 = \norm{u}{L^2_{\alpha-1}(\Omega)}^2  + \norm{\nabla u}{L^2_\alpha(\Omega)}^2$. The space ${V^k_\alpha(\Omega)}$ admits weighted properties analogously to those of the regular Sobolev space, e.g., a weighted version of the Poincaré inequality  \cite{kufner}:
\begin{align*}
\norm{u}{L^2_{\alpha-1}(\Omega)} \leq C_\alpha \norm{\nabla u}{L^2_\alpha(\Omega)} \quad \forall \, u \in V^1_\alpha(\Omega) : u_{\vert_{\partial \Omega}}=0.
\end{align*}

Finally, let us also define the fractional Sobolev space. For $\Omega \subset \mathbb{R}^d$, we define fractional Sobolev space $H^s(\Omega), s\in (0,1)$ using the norm, 
\begin{align*}
\|u\|_{H^s(\Omega)} = \|u\|_{L^2(\Omega)} + \int_\Omega \int_\Omega \dfrac{|u(x) - u(y)|^2}{|x-y|^{d+2s}} \mathrm{d}x \mathrm{d}y.
\end{align*}
We also recall the embedding result of $W^{1,p}(\Omega), \dfrac{2d}{d+2} < p < 2$ in $H^s(\Omega)$ for $s  = 1- d( \frac{1}{p} -\frac{1}{2})$ (see e.g., \cite{adams2003sobolev}, Theorem 7.58). 

\section{Decomposition and regularity properties}

\label{sec:decomposition}
We consider the model in \eqref{eq:basic1}-\eqref{eq:basic2}, posed in a bounded, open domain $\Omega$ with smooth boundary $\partial \Omega$, and Dirichlet boundary conditions;
\begin{subequations}
\begin{align}
 - \nabla \cdot (\kappa \nabla u) &= f \delta_\Lambda  \quad   &\text{in} \, \Omega, \label{eq:model}\\
 u &= u_D  \quad   &\text{on} \, \partial \Omega, \label{eq:model-bc-d}
\end{align}
\end{subequations}
where $\Lambda = \cup_{i=1}^n \Lambda_i$ denotes a collection of line segments $\Lambda_i \subset \Omega$,  $u_D \in C^2(\bar{\Omega})$ some given boundary data, and $\kappa \in W^{2, \infty}(\Omega)$ a uniformly positive, scalar-valued permeability coefficient. In this section we will show that solutions to \eqref{eq:model}-\eqref{eq:model-bc-d} admit a split into singular and high-regularity parts. We start in Section \ref{sec:infinite-line} by illustrating the splitting technique itself and discussing the information it reveals concerning the structure and regularity of the solution. We do this by treating the simple case of $\kappa=1$ and $\Lambda$ consisting of a single line with its endpoints contained in the boundary of $\Omega$. In Section \ref{sec:finite-line}, we then extend the splitting technique to handle arbitrary line sources in the domain, and show that this yields similar results for the regularity. We then proceed to treat the case of $\kappa$ varying in Section \ref{sec:kappa}, before concluding in Section \ref{sec:existence} with a general proof of the splitting properties of solutions to \eqref{eq:basic1}. The proof relies on the discussions in the preceding parts and a simple application of the superposition principle. 
\subsection{The splitting technique}
\label{sec:infinite-line} 
Let us assume $f \in H^2(\Lambda)$, $\kappa=1$ and $\Lambda$ consists of a single line segment passing through the domain $\Omega$, meaning $\Lambda$ has its endpoints contained in $\partial \Omega$. For simplicity, we assume $f=0$ on $\partial \Omega$, meaning that the line source is in-active on the boundary. Next, we orient the coordinate system s.t. $\Lambda$ lies along the $z$-axis, $\Lambda=(0,0,z) \in \Omega$, and employ a cylindrical coordinate system with $r = \sqrt{x^2+y^2}$. The solution to \eqref{eq:model} now admits a splitting into an explicit, low-regularity term $f(z) G(r)$ and a high-regularity correction term $w$:
\begin{align}
u(r,z) = - \frac{1}{2\pi} \left( \underbrace{f(z)  G(r)}_{\text{line singularity}} + \underbrace{w(r,z)}_{\text{correction}} \right),
\label{eq:u-infinite}
\end{align}
where $G = \ln(r)$ and the correction term $w$ is defined as the solution of
\begin{subequations}
\begin{align} 
  -\Delta w &=  f''(z) G(r)     &  \quad \text{in } \Omega \label{eq:w}, \\ 
   w &= -2\pi u_D - f(z) G(r) & \quad \text{on } \partial \Omega, \label{eq:w-bc-d} 
\end{align}
\end{subequations}
The choice of $G=\ln(r)$ is motivated by the observation that $-\ln(r)/2\pi$ is the fundamental solution for $d=2$. Letting $S_{z_0}$ be the intersection of $\Omega$ with the plane perpendicular to $\Lambda$ at a given height $z_0$, i.e., $S_{z_0} = \{ (r, z_0)\in \Omega \} $, it is well known that the Laplacian of $\ln(r)$ in the domain $S_{z_0}$ returns a Dirac point source located in $(0,0,z_0)$, i.e., 
\begin{align}
-\frac{1}{2\pi} \int_{S_{z_0}} \Delta \ln(r)  v \, \mathrm{d} S= v(0,z_0)  \quad  \forall \, v \in C^0(\bar{\Omega}).
\end{align}
It follows then that the Laplacian of $G$ returns the line source on $\Lambda$, i.e.:
\begin{align}
- \frac{1}{2\pi} \int_\Omega \Delta G \, v \, \mathrm{d} \Omega &= - \frac{1}{2\pi} \int_z \int_{S_z} \Delta \ln(r) \, v \, \mathrm{d} S \, \mathrm{d} x \\ &=\int_{\Lambda} v \, \mathrm{d} \Lambda \quad \forall v \in C^0(\bar{\Omega}).
\label{eq:weak-lnr}
\end{align}
The $f(z)G(r)$ term in the decomposition thus acts by prescribing in $u$ a line source of the correct intensity $f$, with the correction function $w$ acting to ensure that $u$ satisfies the given boundary conditions.

The existence of the correction function $w$ follows from standard elliptic theory. By assumption, $f''(z) \in L^2(\Lambda)$, and $G(r) \in L^2(\Omega)$ can be shown by calculation. It follows that the right-hand side in \eqref{eq:w} belongs to $L^2(\Omega)$, and that a $w \in H^2(\Omega)$ exists solving \eqref{eq:w}-\eqref{eq:w-bc-d} in a weak sense. For an example of what each term in the split might look like, the reader is invited to examine Figure \ref{fig:infline-decomposition}. 

To see that $u$ in \eqref{eq:u-infinite} indeed solves \eqref{eq:basic1}, let us calculate its Laplacian:
\begin{align}
- \Delta u &= \frac{1}{2\pi} \left(  \underbrace{f(z) \Delta_\bot G(r)}_{f (z) \delta_\Lambda} +  \underbrace{ 2 \nabla_\perp G(r) \cdot  \nabla_\Vert f}_{=0  \text{ as } \nabla f \perp \nabla G} + G(r) \Delta_\Vert f(z) +  \Delta w(r,z) \right),
\label{eq:laplacian}
\end{align} 
where $\nabla_\Vert$ and $\nabla_\bot$ denote the gradient taken in the line parallel to the line and in the plane perpendicular to the line, respectively; the gradient can then be written as $\nabla = \nabla_\perp + \nabla_\Vert$. By construction, the two last terms in \eqref{eq:laplacian} cancel, and we are left with
\begin{align}
\frac{1}{2\pi}\int_\Omega f(z) \Delta_\bot G(r) \, v \, \mathrm{d}  \Omega = \int_\Lambda f v \, \mathrm{d} \Lambda \quad \forall v \in C^0(\bar{\Omega}).
\end{align}
It follows that the $u$ given by \eqref{eq:u-infinite} solves \eqref{eq:basic1} in a weak sense. 

From the decomposition, we can now investigate the regularity of the solution. Provided that the only source of regularity loss stems from the right-hand side, the correction term $w$ was found to belong to $H^2(\Omega)$; for this reason, we refer to it as being the high-regularity term in the split \eqref{eq:u-infinite}. To see that $f \ln(r)$ is the low regularity term in $u$, let us consider the gradient of $\ln(r)$. A straightforward calculation shows us that $\ln(r) \in L^p(S_{z_0})$ for any $p < \infty$. Similarly, we have $\ln(r) \in L^2_{\alpha-1}$ for any $\alpha>0$. As for the gradient, another calculation shows
\begin{align}
\nabla_\bot \ln(r) =  \frac{(x,y)}{r^2} \in 
\left\{ 
  \begin{aligned}
&  L^{p}(S) & \text{ for } p<2,  \\
&  L^2_\alpha(\Omega) & \text{ for } \alpha>0.
\end{aligned}\right.
\end{align}
In terms of fractional Sobolev space 
we have $W^{1,p}(S), 1 < p < 2$ embedded in $H^s(\Omega)$ for $s  = 1- 2 (\frac{1}{p} -\frac{1}{2})$. Thus, the solution $u$ belonging to $\{ u \in W^{1,p}(\Omega), u_{zz} \in L^2(\Omega)\}$ implies $\{ u \in H^{s}(\Omega),  u_{zz} \in L^2(\Omega)\}$, $s \in (0,1)$. 

To summarize, we have the following results for the regularity of $u$:
\begin{align}
u \in 
\left\{ 
  \begin{aligned}
&  W^{1, p}(\Omega) \cap \lbrace u :  u_{zz} \in L^2(\Omega)  \rbrace & \text{ for } & p<2, \\ 
&  V^1_\alpha(\Omega) & \text{ for } & \alpha>0 \\
&  H^s(\Omega) \cap \lbrace u :  u_{zz} \in L^2(\Omega)  \rbrace & \text{ for } &s \in (0,1).
\end{aligned}\right.
\end{align}
We see that the solution $u$ evades $H^1(\Omega)$, with the regularity loss due to the singular behaviour in the plane normal to $\Lambda$. Moreover, the regularity loss is local to the region surrounding $\Lambda$, meaning that $u$ retrieves full $H^2$ regularity if a small region is removed around $\Lambda$. To be more precise, let $U_\epsilon$ denote a tubular radius surrounding $\Lambda$, 
\begin{align}
U_\epsilon = \lbrace (x,y,z) : \sqrt{x^2+y^2}<\epsilon \rbrace,
\end{align}
and $\Omega_\epsilon$ what remains when this region is removed, i.e., $\Omega_\epsilon = \Omega \setminus U_\epsilon$. We then have $u \in H^2(\Omega_\epsilon)$ for any $\epsilon>0$. Lastly, we see that parallel with $\Lambda$ the solution inherits the regularity of $f$, and therefore exhibits high regularity in the $z$-direction.

\subsection{Decomposition and regularity for arbitrary line sources}
\label{sec:finite-line}
The decomposition presented in the previous section is valid for the case $f \in H^2(\Lambda)$ with the line $\Lambda$ passing through the domain and $f=0$ on $\Lambda \cap \partial \Omega$. This excludes, however, a number of interesting scenarios, e.g., when $\Lambda$ is allowed to branch. For this reason, let us now extend the decomposition technique to handle $\Lambda$ being some arbitrary line segment between the points $\textbf{a}, \textbf{b} \in \Omega$, keeping the assumptions of $\kappa=1$ and $f\in H^2(\Lambda)$. 

The line $\Lambda$ can be described by the parametrization $\mathbf{y} = \mathbf{a} + \pmb{\tau} t \quad \text{for } t  \in (0, L),$
where $L = \Vert \mathbf{b}-\mathbf{a} \Vert$ denotes the Euclidean norm and $\pmb{\tau} = (\mathbf{b} - \mathbf{a}) / L$ is the normalized tangent vector of $\Lambda$. As in Section \ref{sec:infinite-line}, we want to construct a function $G$ so that its Laplacian satisfies
\begin{align}
\int_\Omega - \Delta G \, v \, \mathrm{d} \Omega = \int_{\Lambda} v \, \mathrm{d} \Lambda \quad \forall v \in C^0(\bar{\Omega}).
\end{align} Considering for the moment $\Omega = \mathbb{R}^3$, the Green's function for the Laplacian is given by
\begin{align}
G_{3D}(\mathbf{x}, \mathbf{y}) = \frac{1}{4\pi} \frac{1}{\Vert \mathbf{x}-\mathbf{y} \Vert},
\end{align}
and a candidate $G$ can now be found by the convolution of $\delta_\Lambda$ and $G_{3D}(\mathbf{x}, \mathbf{y})$:
\begin{align}
G(\mathbf{x}) &=  \int_\Omega \frac{\delta_{\Lambda}}{ \Vert \mathbf{x} - \mathbf{y} \Vert} \, \mathrm{d} \mathbf{y} \\
&=\int_0^L \frac{1}{ \Vert \mathbf{x} - (\mathbf{a} + \pmb{\tau} t) \Vert} \, \mathrm{d}t  \\
&= \ln \left(    \frac{r_{b} + L +  \pmb{\tau} \cdot (\mathbf{a}-\mathbf{x})   } {r_{a} + \pmb{\tau} \cdot (\mathbf{a}-\mathbf{x})   }  \right),
\label{eq:G}
\end{align}
where $r_{b}(\mathbf{x}) =  \Vert \mathbf{x} - \mathbf{b} \Vert$ and $r_{a}(\mathbf{x})  =\Vert \mathbf{x} - \mathbf{a} \Vert$.

Solving now the line source problem with $f$ as the line source intensity, the solution admits the decomposition
\begin{align}
u = \frac{1}{4\pi} \left( \underbrace{ E(f) G }_{\text{line singularity}} + \underbrace{w}_{\text{correction}} \right),
\label{eq:finite-line-decomp}
\end{align}
where $E: H^2(\Lambda) \rightarrow H^2(\Omega) \cap C^2(\bar{\Omega})$ denotes a suitable extension of the line source intensity $f$ into the domain $\Omega$. Depending on the choice of extension, the correction term $w$ must solve
\begin{subequations}
\begin{align} 
  - \Delta w &= F    &  \quad \text{in } \Omega,    \label{eq:wfinite}\\ 
   w &= 4 \pi u_D - f G & \quad \text{on } \partial \Omega., 
\end{align}
\end{subequations}
with right-hand side $F$ given by
\begin{align}
F = \Delta E(f) G + 2 \nabla E(f) \cdot \nabla G.
\end{align}
The $1/4\pi$-scaling of $G$ in \eqref{eq:finite-line-decomp} differs from the $-1/2\pi$-scaling in Section \ref{sec:infinite-line} as we now used the fundamental solution for $d=3$. To be more precise, the scalings $1/2\pi$
and $1/4\pi$ enter in the fundamental solution as the volume of a unit ball in $\mathbb{R}^d$ for dimensions $d=2$ and $d=3$, respectively (see e.g. \cite[p. 22]{evans10}). The relation
\begin{align*}
\lim_{L \rightarrow \infty} \, \frac{1}{4\pi}   \ln \left(    \frac{r_{b} + L +  \pmb{\tau} \cdot (\mathbf{a}-\mathbf{x})   } {r_{a} + \pmb{\tau} \cdot (\mathbf{a}-\mathbf{x})   }  \right) \approx - \frac{1}{2\pi} \ln(r);
\end{align*}
is used in electromagnetism to approximate the potential of an infinite length line charge. The splitting with $-\ln(r)/2\pi$ introduced in Section \ref{sec:infinite-line} can thus be thought of to hold for any infinite length line segment. It is further valid for a smooth line source intensity $f$ as $f$ is then $H^2$-extendible to the infinite extension of the line segment. The same does not hold if $f$ is discontinuous along $\Lambda$; for this case, it is necessary to use $G$ as it is given in \eqref{eq:G}.

\begin{figure}
\centering
\includegraphics[width=0.18\textwidth]{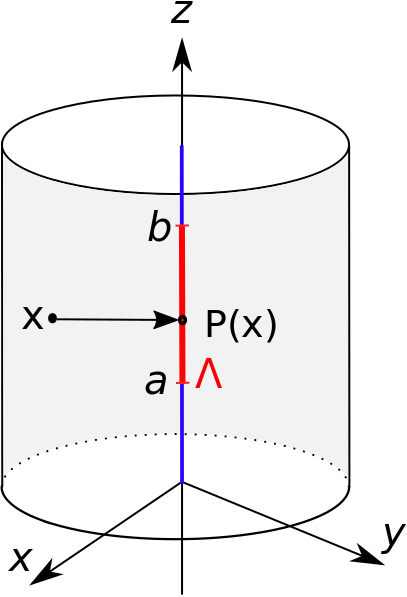}
\caption{A cylindrical domain $\Omega$ and the line segment $\Lambda=\{ (0,0,z) : z \in [a,b]\}$, its elongation $\Lambda= \lbrace (0,0,z) \in \Omega \rbrace$, and a projection operator $P:\Omega \rightarrow \ \tilde{\Lambda}$ taking a point in the domain to the closest point on $\tilde{\Lambda}$.}
\label{fig:finite-line}
\end{figure}

The regularity of the correction function $w$ thus depends now on the choice of extension $E$; for $w \in H^2(\Omega)$, it is necessary to require $F \in L^2(\Omega)$. By a calculation, one finds that $\nabla G_\bot$ fails to be $L^2$-integrable in the neighbourhood of $\Lambda$. The result $w \in H^2(\Omega)$ thus requires choosing $E$ so that the term involving $\nabla G_\bot$ is cancelled; this can be achieved, for example, by choosing $E(f)$ s.t. $\nabla_\perp E(f) = (0,0)$ in $U_\epsilon$. Assuming $f=f(s)$ is a given, analytic function with respect to the arc length parameter $s$, the extension operator can be chosen as
\begin{align}
E(f)(\mathbf{x}) =  f(P(\mathbf{x})  ) & \text{ for } \mathbf{x} \in \Omega
\label{eq:extension}
\end{align}
where $P : \Omega \rightarrow \tilde{\Lambda}, \, \mathbf{x} \rightarrow (\mathbf{x}-\mathbf{a}) \cdot \pmb{\tau}$ denotes the orthogonal projection of a point $\mathbf{x}$ onto the elongation of $\Lambda$, as is illustrated in Figure \ref{fig:finite-line}.
For an example of what the splitting might look like with this choice of extension, the reader is invited to examine Figure \ref{fig:finite-line-decomposition}. 

Given the extension \eqref{eq:extension}, the right-hand side $F$ can be stated more explicitly. To see this, let us first note that
\begin{align}
\nabla G &= \frac{\frac{\mathbf{x}-\mathbf{b}}{r_b} - \pmb{\tau}  }{L + r_b + \pmb{\tau}\cdot (\mathbf{a}-\mathbf{x})} -  \frac{\frac{\mathbf{x}-\mathbf{a}}{r_a} - \pmb{\tau}  }{r_a + \pmb{\tau}\cdot (\mathbf{a}-\mathbf{x})},
\end{align}
and that the gradient and Laplacian of $E(f)$ are given by $\nabla E(f) = f'(P(\mathbf{x}))  \pmb{\tau} $ and $\Delta E(f) = f''(P(\mathbf{x}))  \Vert \pmb{\tau} \Vert = f''(P(\mathbf{x}))$. It follows that $F$, for this particular choice of extension operator, is given by
\begin{align}
F = f''(p(\mathbf{x})) G + f'(p(\mathbf{x})) \left(\frac{1}{r_a} - \frac{1}{r_b} \right).
\label{eq:w2}
\end{align}

To see that $F \in L^2(\Omega)$, let us for simplicity assume the domain to be cylindrical around $\Lambda$, and shift the coordinate system so that $\Lambda$ lies along the $z$-axis, $\Lambda= \lbrace (0,0,z): \, \, z \in [a, b]  \rbrace.$ The extension $E(f)$ can then be expressed more simply as $E(f) = f(z)$ and the function $G$ as
\begin{align}
G(r, z) 
 &  \frac{1}{4\pi} \ln \left( \frac{  r_b-(z-b)} { r_a-(z-a) }\right) \quad \text{ where } r_a = \sqrt{r^2+(z-a)^2}  \quad \text{ and }  r_b = \sqrt{r^2+(z-b)^2}.
\end{align}
Without loss of generality, let us shift the coordinate system so that the endpoint $(0,0,b)$ is the origin: a calculation then shows us that $f'(z)/r_b \in L^2(\Omega)$:
\begin{align}
\norm{f'(z) / r_b}{L^2(\Omega)}^2 &= \int_0^{2\pi}\int_0^{Z} f'(z)^2  \int_0^{R} \frac{1}{r^2+z^2}  r \,  \mathrm{d} r  \mathrm{d} z \, \mathrm{d} \theta   \\
&= 2\pi \int_0^Z f'(z)^2 \ln\left( \frac{R^2}{z^2} + 1\right) \, \mathrm{d} z < \infty,
\end{align}
and that $f(z) \ln(r_b-(z-b)) \in L^2(\Omega)$:
\begin{align}
\norm{\ln(r_b - z)}{L^2(\Omega)}^2 &= \int_0^{2\pi}\int_0^{Z} f(z)^2  \int_0^{R}  \ln(\sqrt{r^2+z^2}-z)^2r \,  \mathrm{d} r  \mathrm{d} z  \mathrm{d} \theta < \infty. 
\end{align}
The terms depending on $r_a$ can be similarly treated by shifting the origin to $(0,0,a)$. 

It only remains to understand the regularity properties of the solution. As before, the correction term $w$ and extension $E(f)$ both belong to $H^2(\Omega)$, and the low regularity terms in $u$ are now $\ln(r_b-(z-b))$ and $\ln(r_a-(z-a))$. Both of these logarithmic terms belong to $L^2_{\alpha-1}(\Omega)$. To understand the regularity of their gradients, considering the domain $\Omega / \Lambda$, in which case the gradient of $\ln(r_b-(z-b))$ is given by
\begin{align}
\nabla \ln(r_b-z) = \left( \frac{x}{r^2+z^2-z\sqrt{r^z+z^2}}, \frac{y}{r^2+z^2-z\sqrt{r^z+z^2}}, -\frac{1}{r_b} \right).
\end{align}
Shifting again the coordinate system so that $(0,0,b)$ is the origin, we can examine the norm of this gradient in a weighted $L^2_\alpha(\Omega)$. A straightforward calculation reveals
\begin{align}
\norm{\nabla_\Vert \ln(r_b-z)}{L^2_\alpha(\Omega)}^2 &= \int_\Omega \frac{r^{2+\alpha}}{(r^2+z^2-z\sqrt{r^z+z^2})^2} \, \mathrm{d}  \Omega \\
&= \int_0^{2\pi} \int_0^{R} \int_a^{b} \frac{r^{3+\alpha}}{(r^2+z^2-z\sqrt{r^z+z^2})^2} \, \mathrm{d} z \, \mathrm{d} r \, \mathrm{d}\theta  \\
&={2\pi} \int_0^{R}    
\left[ 2 r^{\alpha-1} \left(\sqrt{z^2+r^2}+z\right)- r^\alpha\arctan\left(\dfrac{z}{r}\right) \right]_a^0 
\mathrm{d} r \\
&= {2\pi} \int_0^{R}    2 r^{\alpha-1} \left(r - \sqrt{a^2+r^2}-a\right)+ r^\alpha\arctan\left(\dfrac{a}{r}\right)
\mathrm{d} r  \\
&< \infty \quad \text{ if } \alpha > 0.
\end{align}
A similar argument holds for $\nabla \ln(r_a-(z-a))$, and for the $L^p$-norm with $p<2$. For the $z$-component of these gradients, recall $1/r_b, 1/r_a \in L^2(\Omega)$. Notice however, that $u_{zz}$ falls just short of belonging to $L^2(\Omega)$. The regularity properties are therefore as follows:

\begin{align}
u \in
\left\{ 
  \begin{aligned}
&  W^{1, p}(\Omega) \cap \left \{ u_{zz} \in L^{p}(\Omega) \right \} & \text{ for } & p<2,  \\
&  V^1_\alpha(\Omega) & \text{ for } & \alpha>0, \\
& H^s(\Omega) \cap \left \{ u_{zz} \in L^{p}(\Omega) \right \} & \text{ for } & s<1, p<2.
\end{aligned}\right.
\end{align}

As in the previous section, the solution $u$ thus barely eludes $H^1(\Omega)$. The regularity loss is again local to the region surrounding $\Lambda$, meaning $u \in H^2(\Omega_\epsilon)$ for any $\epsilon>0$. Moreover, the solution exhibits anisotropic regularity properties, with improved regularity parallel to the line.

\subsection{Spatially varying permeability}
\label{sec:kappa}
Let us now extend the splitting technique to handle spatially varying $\kappa$. Assuming $\kappa \in W^{2,\infty}(\Omega)$ and $\kappa >0$, the solution $u$ to \eqref{eq:basic1} admits the splitting
\begin{align}
u = \frac{1}{4\pi} \left( \underbrace{ \frac{ E(f) G}{\kappa}}_{\text{line singularity}} + \underbrace{w}_{\text{correction}} \right),
\label{eq:u-infinite-line-kappa}
\end{align}
with the correction term $w$ solving
\begin{subequations}
\begin{align} 
  - \nabla \cdot (\kappa \nabla w) &= F    &  \quad \text{in } \Omega, \label{eq:w-convdiff1} \\ 
   w &= 4\pi u_D - f G & \quad \text{on } \partial \Omega \label{eq:w-convdiff1-bc-d}
\end{align}
\end{subequations}
with right-hand side
\begin{align}
F =
\left( \Delta E(f) - \frac{ E\left( f \right) \Delta \kappa + \nabla E(f) \nabla \kappa }{\kappa}  + \frac{ E(f) \vert \nabla\kappa \vert^2}{\kappa^2} \right)  G +  \left( 2 \nabla E(f)-  \frac{E(f) \nabla \kappa}{\kappa^2} \right) \cdot & \nabla G 
\end{align}
To see that \eqref{eq:u-infinite-line-kappa} solves \eqref{eq:basic1}, we calculate $-\nabla \cdot (\kappa \nabla u)$. By construction, this yields the single term
\begin{align}
-\nabla \cdot (\kappa \nabla u) = - \frac{1}{4\pi} f \Delta G.
\end{align}
Furthermore, we know from Section \ref{sec:finite-line} that this term returns the line source with required intensity $f$.

The regularity of $w$ depends now on the behaviour of $\kappa$. Recall that $\nabla_\perp G$ falls just short of belonging to $L^2(\Omega)$, and consequently
\begin{align}
F \in
\left\{ 
  \begin{aligned}
&  L^{p}(\Omega)   & \text{ for } p<2,  \\
&  L^2_\alpha(\Omega) & \text{ for } \alpha>0 ,
\end{aligned}\right. \qquad \Rightarrow \qquad
w \in
\left\{ 
  \begin{aligned}
&  W^{2, p}(\Omega)  & \text{ for } p<2,  \\
&  V^1_\alpha(\Omega) & \text{ for } \alpha>0,\\
&  H^k(\Omega) & \text{ for } k<2. 
\end{aligned}\right.
\end{align}

We see that for general $\kappa \in W^{2,p}(\Omega)$ the correction function $w$ falls just short of belonging to $H^2(\Omega)$. Notice, however, that it still constitutes a higher regularity term compared to the full solution $u \, \cancel{\in} \, H^1(\Omega)$. Moreover, given the existence of some $\epsilon>0$ such that $\nabla_\perp \kappa=0$ in the small tubular neighbourhood $U_\epsilon$, $R$ will indeed belong to $L^2(\Omega)$ and we recover $w \in H^2(\Omega)$.

\subsection{Splitting theorem}
\label{sec:existence}
We are now ready to summarize the results of Sections \ref{sec:infinite-line}-\ref{sec:kappa} with a splitting theorem. The results hold for any collection of line sources $\delta_{\Lambda}$ when $\Lambda$ is a collection of line segments and $f$ is piecewise $H^2$ on each line segment. The line source intensity $f$ is thus allowed to contain jumps. The jump can be handled by splitting the line segment containing it into two pieces. By superposition, these two line segments can be handled separately  using the splitting technique shown in Section \ref{sec:finite-line}.

\begin{thrm}[Singularity splitting theorem for elliptic equations with line sources]
Let $\Lambda= \cup_{i=1}^n\Lambda_i$ be a collection of line segments $\Lambda_i$, $f$ be a piecewise $H^2$ function on each line segment $\Lambda_i$, and $\kappa \in W^{2,p}(\Omega)$ be a positive, scalar-valued permeability. The solution $u \in W^{1,p}(\Omega)$ for $p<2$ solving \eqref{eq:basic1} then admits the split
\begin{align}
u(\mathbf{x}) = \frac{1}{4\pi} \left( \sum_{i=1}^n \frac{  E_i(f)  G_i}{\kappa}  + w  \right),
\label{eq:u-existence}
\end{align}
where $G_i$ is the logarithmic term
\begin{align}
G_i(\mathbf{x})  = \ln \left(    \frac{r_{b, i} + L_i +  \pmb{\tau}_i \cdot (\mathbf{a}_i-\mathbf{x})   } {r_{a, i} + \pmb{\tau}_i \cdot (\mathbf{a}_i-\mathbf{x})   }  \right),
\end{align}
$E_i(f)$ is an extension operator $E_i(f): H^2(\Lambda) \rightarrow H^2(\Omega)$ that extends $f$ from line segment $\Lambda_i$ into the domain, assumed to satisfy 
\begin{align*}
(i) \quad& E_i(f)(\mathbf{x}) = f(\mathbf{x}) \text{ for } \mathbf{x} \in \Lambda_i, \\
(ii) \quad& E_i(f) \in H^2(\Omega), \\
(iii) \quad& \nabla G_i \cdot \nabla E_i(f) \in L^2(\Omega),
\end{align*}
and $w$ solves 
\begin{align} 
   - \nabla \cdot ( \kappa\nabla w) &= F    &  \quad \text{in } \Omega, \\ 
   w &= w_D & \quad \text{on } \partial \Omega,
\end{align}
with right-hand side $F$ given by
\begin{align}
F & =  \sum_{i=1}^n   \left( \Delta E_i(f) -  \frac{E_i(f) \Delta \kappa +  \nabla E_i(f) \cdot \nabla \kappa }{\kappa}  + \frac{E_i(f) \vert \nabla\kappa \vert^2}{\kappa^2} \right) G_i  +  \left( 2 \nabla E_i(f) -  \frac{ E_i(f) \nabla \kappa}{\kappa}  \right) \cdot \nabla G_i, 
\label{eq:R}
\end{align}
and boundary data $w_D$ given by
\begin{align}
w_D &= 4\pi u_D - \sum_{i=1}^n  E_i(f) G_i. \label{eq:bcs1}
\end{align}

\label{thm:existence}

\end{thrm}

\begin{proof}
A direct calculation shows that, for $u$ being the solution given in \eqref{eq:u-existence},
\begin{align}
-\nabla \cdot (\kappa \nabla u) = \frac{1}{4\pi} \Bigg( \sum_{i=1}^n  E_i(f) \Delta G_i + \Delta E_i(f) G_i + 2 \nabla E_i(f) \cdot \nabla G_i - \nabla \cdot \left( \frac{\nabla \kappa}{\kappa}E_i(f) G_i \right) -\nabla \cdot (\kappa \nabla w) \Bigg)
\end{align}
By construction, all these terms cancel except the first, and we have
\begin{align}
-\nabla \cdot (\kappa \nabla u) &= \frac{1}{4\pi} \left( \sum_{i=1}^n  E_i(f) \Delta G_i ) \right).
\end{align}
Recalling now from Section \ref{sec:finite-line} that
\begin{align}
\int_\Omega -(\Delta G_i) v \, \mathrm{d} \Omega = \int_{\Lambda
_i} v \, \mathrm{d} \Lambda_i \quad \forall v \in C^0(\bar{\Omega}),
\end{align}
we see that
\begin{align}
-\nabla \cdot (\kappa \nabla u) &= f \delta_\Lambda,
\end{align}
where the Dirac distribution on the line is understood in the sense of \eqref{eq:basic2}. Moreover, by construction, $u$ matches the prescribed boundary conditions. It follows that the $u$ constructed in \eqref{eq:u-existence} solves \eqref{eq:model}-\eqref{eq:model-bc-d} in a weak sense. 
\end{proof}

\section{Numerical methods}

\label{sec:numerical-disc}
Let us now consider the numerical approximation of \eqref{eq:model}-\eqref{eq:model-bc-d} by finite element methods. From this point on, we consider only $\kappa=1$. For the numerical discretization we will make use of both simplicial and prismatic Lagrange elements, the latter as it lets us study the anisotropic properties of the solution. We consider two different numerical methods, where the first is the standard lowest-order Galerkin FE method \eqref{eq:model} introduced in \cite{dangelo2012}. The second method is based on the splitting properties of the solution: Using Theorem \ref{thm:existence}, we can formulate a numerical method in which $u$ is approximated via its correction function $w$. We shall refer to this method as the Singularity Subtraction Based Finite Element (SSB-FE) method for the line source problem.

\subsection{Discretization}
Let us assume the domain $\Omega$ is a polyhedron that readily admits a partitioning $\mathcal{T}_{K, h}$ into simplicials $K$,
\begin{align*}
\bar{\Omega} = \bigcup_{K \in \mathcal{T}_{K, h}} K,
\end{align*}
as well as a partitioning $\mathcal{T}_{P, h}$ into prisms $P$,
\begin{align*}
\bar{\Omega} = \bigcup_{P \in \mathcal{T}_{P, h}} P.
\end{align*}
The prismatic elements $P \in \mathcal{T}_{P,h}$ consist of the  product $T \times I$ between a triangle $T$ in the $x,y$-plane and an interval $I$ in the $z$-axis. The simplicial mesh is characterized by the mesh discretization parameter $h$, taken as maximum element size $h = \max_{K \in \mathcal{T}_{K,h}} h_K$. The prismatic element is characterized by the two mesh  discretization parameters $h_\perp$ and $h_\Vert$, where $h_\perp = \max_{T \times I  \in \mathcal{T}_P,h} h_{T}$ measures the size of the triangles in the $x,y$-plane, and $h_\Vert = \max_{T \times I \in \mathcal{T}_{P,h}} h_{I}$ measures the discretization of the interval along the $z$-axis. Both meshes $\mathcal{T}_{K, h}$ and  $\mathcal{T}_{P, h}$ are assumed to satisfy all the requirements of a conforming mesh.

The partitionings can now be associated with discrete spaces. For the simplicial partitioning $\mathcal{T}_{K, h}$ we pick the standard conforming space
\begin{align*}
V^{K, h}_{u_D} = \lbrace v_h \in C^0_{u_D}(\Omega), \, v_h \vert_K \in \mathbb{P}_1  \text{ where } K \in \mathcal{T}_{K, h} \rbrace,
\end{align*}
where $\mathbb{P}^1$ denotes the space of polynomials of degree $1$ and $C^0_{u_D}(\Omega)$ the space of continuous elements that equal the interpolation of $u_D$ on the boundary, i.e., 
\begin{align*}
C^0_{u_D}(\Omega) = \lbrace u \in C^0(\Omega) : u\vert_{\partial \Omega} = \mathcal{I}_h u_D \rbrace.
\end{align*} 
Similarly, for the prismatic partitioning $\mathcal{T}_{P,h}$ we pick the conforming space
\begin{multline}
V^{P, h}_{u_D} = \lbrace v_h \in C^0_{u_D}(\Omega), \, v_h(x,y,z) \vert_P = \sum_{i=1}^3 \sum_{j=2}^3  \sigma_{ij} \lambda_i(x,y) l_j(z),  \\
\text{ where } P \in \mathcal{T}_{P, h}, P = T \times I, \sigma_{ij} \in \mathbb{R}, \lambda_i \in\mathbb{P}^1(T), l_j \in \mathbb{P}^1(I)   \rbrace,
\end{multline}
defined on each prism $P = T \times I$ as the Cartesian product between linear elements defined on the triangle $T$ and linear elements defined on each interval $I$. Similar to the simplicial elements, the degree of freedoms for the prismatic element are defined on its vertices, as is illustrated in Figure \ref{fig:prism}.

The first order Galerkin approximation of \eqref{eq:basic1}, which we shall denote as the standard FE method, can now be introduced: Find $u_h \in V^h_{u_D}$ s.t.
\begin{align}
(\nabla u_h, \nabla v_h)_\Omega &=  (f, v_h)_\Lambda  &  \label{eq:varform-uh} \text{ for all }  v_h \in V^h_0,
\end{align}
where $V^h$ can be either $V^{K,h}$ or $V^{P,h}$. The stability of this approximation was shown in \cite{dangelo2012} for the simplicial type space $V^{K,h}$, and we choose here to assume this can be extended to hold also for the prismatic type space $V^{P,h}$.

The line source in the right-hand side of \eqref{eq:varform-uh} is well known to cause sub-optimal convergence rates for the approximation. For this reason, we define the SSB-FE method, in which $u$ is solved for via its correction function $w$: 
Find $w_h \in V^h_{w_D}$ such that
\begin{align}
(\nabla w_h, \nabla v_h)_\Omega &= (F, v_h)_\Omega &  \label{eq:varform-wh} \text{ for all }  v_h \in V^h_{0},
\end{align}
where the right-hand side $F$ and boundary data $w_D$ are given by \eqref{eq:R} and \eqref{eq:bcs1}, respectively. 

\subsection{Error estimates}
The logarithmic term in the analytic solution of $u$ is well known to cause sub-optimal convergence rates when solving $u_h$ directly via \eqref{eq:varform-uh}. In \cite{Scott1973}, it was found that approximating $\ln(r)$ near the origin with piecewise linears yields errors in the $L^2(\Omega)$ norm of order $\mathcal{O}(h^{1-\epsilon})$, $\epsilon>0$ being some arbitrarily small parameter. In \cite{dangelo2012}, the following error estimate was found for the approximation \eqref{eq:varform-uh} with simplicial elements:
\begin{align*}
\norm{u-u_h}{V^1_\alpha(\Omega)} \leq C_\alpha h^{\alpha - \epsilon} \norm{u}{V^2_{\alpha+1}(\Omega)},
\end{align*} 
where $\alpha$ is assumed to satisfy $0<\alpha<1$ and $\epsilon \in (0, \alpha)$. Note that convergence is not possible in the standard $H^1$-norm as $u$ fails to belong to $H^1(\Omega)$. The Aubin-Nitsche theorem in weighted norms then predicts convergence of order $1+\alpha-\epsilon$ in the $L^2_\alpha(\Omega)$ norm; in particular, we expect the sub-optimal error estimate 
\begin{align*}
\norm{u-u_h}{L^2(\Omega)} \lesssim h^{1-\epsilon} \norm{u}{H^1_\alpha(\Omega)}.
\end{align*} 
in the standard $L^2$-norm, where $X \lesssim Y$ is taken to denote $X \leq C Y$ for some constant $C>0$.

As $w$ belongs to $H^2(\Omega)$, the approximation of $w$ by \eqref{eq:varform-wh} will yield optimal convergence rates:
\begin{align}
\norm{w-w_h}{L^2(\Omega)} &\lesssim h^{2} \norm{w}{H^2(\Omega)}, \\
\norm{w-w_h}{H^1(\Omega)} &\lesssim h^{1} \norm{w}{H^2(\Omega)}.
\end{align}

\subsection{Modelling error}
Let us now consider the effect of removing a collection of line segments from $\Lambda$. The motivation for such a removal may be to improve simulation runtime by removing line segments that offer a negligible contribution to the total solution. Another motivation might be to assess the effect of imperfect data acquisition, or to investigate the effect a disease such as stroke has in altering the blood flow to the surrounding tissue. 

We start by defining $u_0$ to be the solution of the boundary-value problem
\begin{align}
- \Delta u_0 &= 0 \quad   &\text{in} \, \Omega, \\
 u_0 &= u_D  \quad   &\text{on} \, \partial \Omega,
\end{align}
and $u_i$, $i=1,2,...,n$, to be the solution found when considering line source $i$, i.e.,
\begin{align}
- \Delta u_i &= f \delta_{\Lambda_i} \quad   &\text{in} \, \Omega, \label{eq:model-red}\\
 u_i &= 0  \quad   &\text{on} \, \partial \Omega. \label{eq:model-bc-red}
\end{align}
The total solution $u$, solving the full system \eqref{eq:basic1}-\eqref{eq:basic2}, can then be found by summation:
\begin{align*}
u_{n} = u_0 + \sum_{i=1}^{n} E_i(f) G_i + w_i,
\end{align*}
where $w_0$ solves the boundary value problem 
\begin{align}
- \Delta w_0 &= 0 \quad &\text{in} \, \Omega,\\
 w_0 &= w_D \quad   &\text{on} \, \partial \Omega, 
\end{align}
with $w_D = - 4\pi u_D$. The remaining terms $w_i$, $i=1,...,n$, each constitute the correction term associated with line source $i$, and are defined as the solution to  
\begin{align}
- \Delta w_i &= F_i \quad &\text{in} \, \Omega, \label{eq:model-red-w}\\
 w_i &= - E_i(f) G_i \quad   &\text{on} \, \partial \Omega, \label{eq:model-bc-red-w}
\end{align}
with 
\begin{align}
F_i = \Delta E_i(f) G_i + 2 \nabla E_i(f) \cdot \nabla G_i.
\end{align}

Next, we define $u_{n_0}$ to be the solution found when considering only $n_0-1 < n = n_a + n_v$ of the line sources:
\begin{align*}
u_{n_0} = \frac{1}{4\pi} \left( \sum_{i=0}^{n_0-1} E_i(f) G_i + w_i \right).
\end{align*}
The error $\norm{u-u_h}{L^2(\Omega)}$ will then consist of two terms,
\begin{align*}
\norm{u-u_h}{L^2(\Omega)} & \leq \norm{u_{n_0}-u_h}{L^2(\Omega)} +  \norm{u-u_{n_0}}{L^2(\Omega)} ,
\end{align*}
where the first term is the usual numerical error, with $u_{h}$ being the FE approximation of $u_{n_0}$, and the second term is the modelling error. From the splitting properties of the solution, we find that this modelling error satisfies
\begin{align}
\norm{u-u_{n_0}}{L^2(\Omega)} & \lesssim \sum_{i={n_0}}^{{n}} \norm{E_i(f) G_i + w_i }{L^2(\Omega)} 
\\
& \lesssim \sum_{i={n_0}}^{{n}} \norm{E_i(f) G_i }{L^2(\Omega)} +  \norm{F_i}{L^2(\Omega)}  +\norm{w_i}{H^{\frac{1}{2}}(\partial \Omega)} \\
& \lesssim  \sum_{i={n_0}}^{{n}} \norm{E_i(f) G_i }{L^2(\Omega)} + \norm{E_i(f) G_i }{H^\frac{1}{2}(\partial \Omega)} + \norm{\Delta E_i(f) G_i}{L^2(\Omega)} + 2 \norm{\nabla E_i(f) \cdot \nabla G_i}{L^2(\Omega)}
\label{eq:model-error}
\end{align}
While this estimate is not sharp, it suggests that the modelling error introduced when neglecting line segment $i$ depends on two factors: the line source intensity $f_i$, as well as the logarithmic term $G_i$ associated with it. The exact interplay between the two is hard to quantify; care must therefore be taken when considering the impact of removing a line segment. We will return to this thought in Section \ref{sec:num-brain}.

\section{Numerical results}
The purpose of this section is twofold. Firstly, we want to verify numerically the regularity properties and error estimates found in Sections \ref{sec:decomposition} and \ref{sec:numerical-disc}, respectively. To this end, we have performed numerical experiments on the simple unit cube domain $\Omega = (0, 1)^3$, with $\Lambda$ taken as a single, vertical line. Physically, these experiments can be interpreted as modelling a single well injecting fluid into a box-shaped reservoir. We test the two different numerical approaches presented in Section \ref{sec:numerical-disc}, the first one being the straightforward FE method described by equations \eqref{eq:varform-uh}, and the second being the SSB-FE method given by \eqref{eq:varform-wh}. Furthermore, we utilize here a prismatic mesh; this allows us to refine independently in the directions parallel with and perpendicular to the line. Prismatic meshes are more commonly used in reservoir engineering as it allows for a finer discretization in the horizontal plane; in particular, it allows for easy grading of the mesh around a well \cite{Mundal2010}. Here, we will use it to probe the anisotropic properties of the solutions.  The computations for this first part were performed using the finite element framework Firedrake \cite{Dalcin2011}, \cite{Rathgeber2016}, which relies on PetSc \cite{petsc-user-ref}, \cite{petsc-efficient}. 

The second purpose of this section is to demonstrate the capabilities of the SSB-FE method in handling problems with a large number of line segments. To this end, we test it on a dataset describing the vascular system in a human brain. As this dataset contains nearly 3000 line segments, it would be computationally challenging to resolve using the standard FE method. We show that our numerical approach for this test case is advantageous, in that it allows for a fine resolution of the pressure profile around each line segment and facilitates model reduction techniques. Here, the computations were performed using the finite element framework FEniCS \cite{LoggMardalEtAl2012a}.

\begin{figure}
\centering
\hspace{4em}
\includegraphics[width=0.18\textwidth]{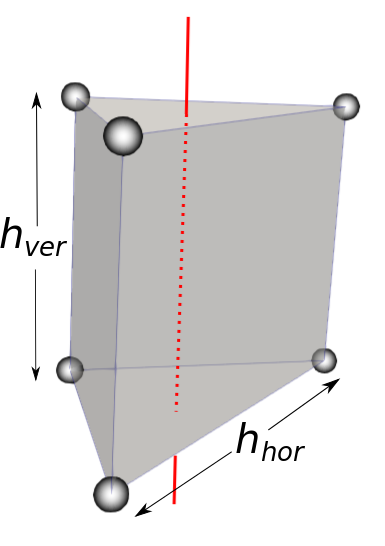}
\includegraphics[width=0.34\textwidth]{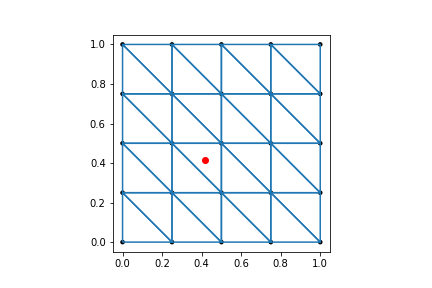}
\caption{(right) The prismatic element used in the finite element methods, with its degrees of freedom indicated by grey spheres, and line placement in red. (left) Line placement in a horizontal cross-section of the coarsest mesh.}
\label{fig:prism}
\end{figure}

\subsection{Convergence test for smooth f}
\label{sec:num-infinite-line}

\renewcommand{\thetable}{\arabic{table-1}} 

\begin{figure*}
    \centering
    \begin{subfigure}[t]{0.34\textwidth}
        \centering
        \includegraphics[height=1.45in]{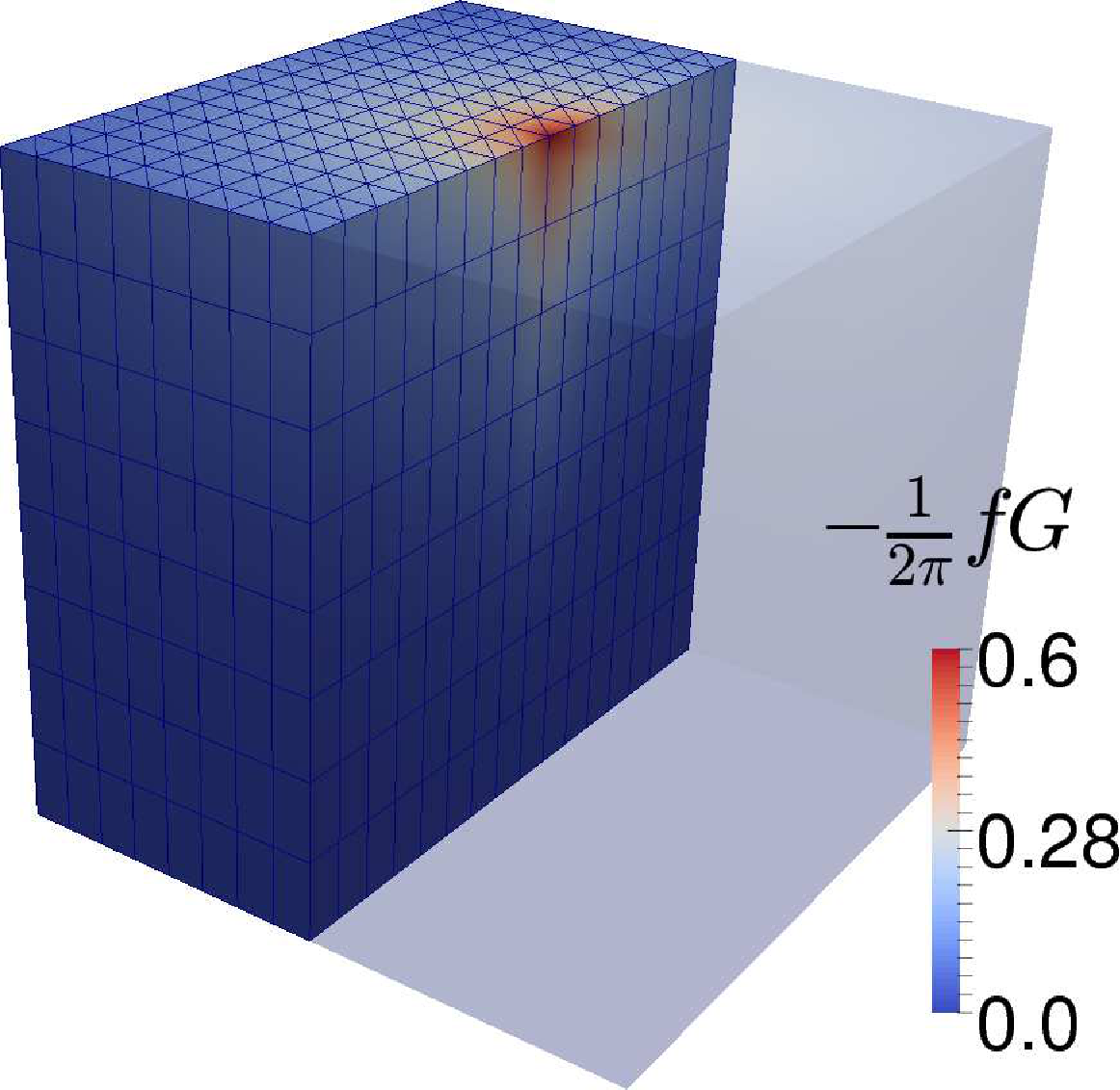}
        \caption{ $ - \frac{1}{2\pi} f(z) \ln(r)$}
    \end{subfigure}
    \begin{subfigure}[t]{0.32\textwidth}
        \centering
        \includegraphics[height=1.45in]{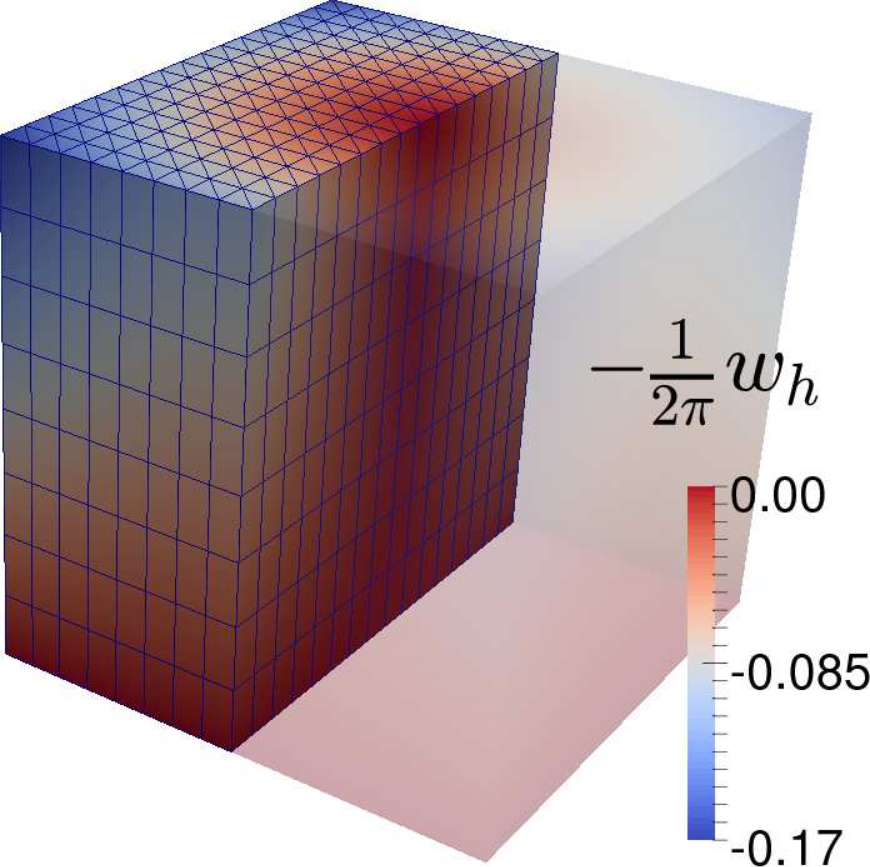}
        \caption{ Correction term: $ -\frac{1}{2\pi} w_h $}
    \end{subfigure}
    \begin{subfigure}[t]{0.32\textwidth}
        \centering
        \includegraphics[height=1.45in]{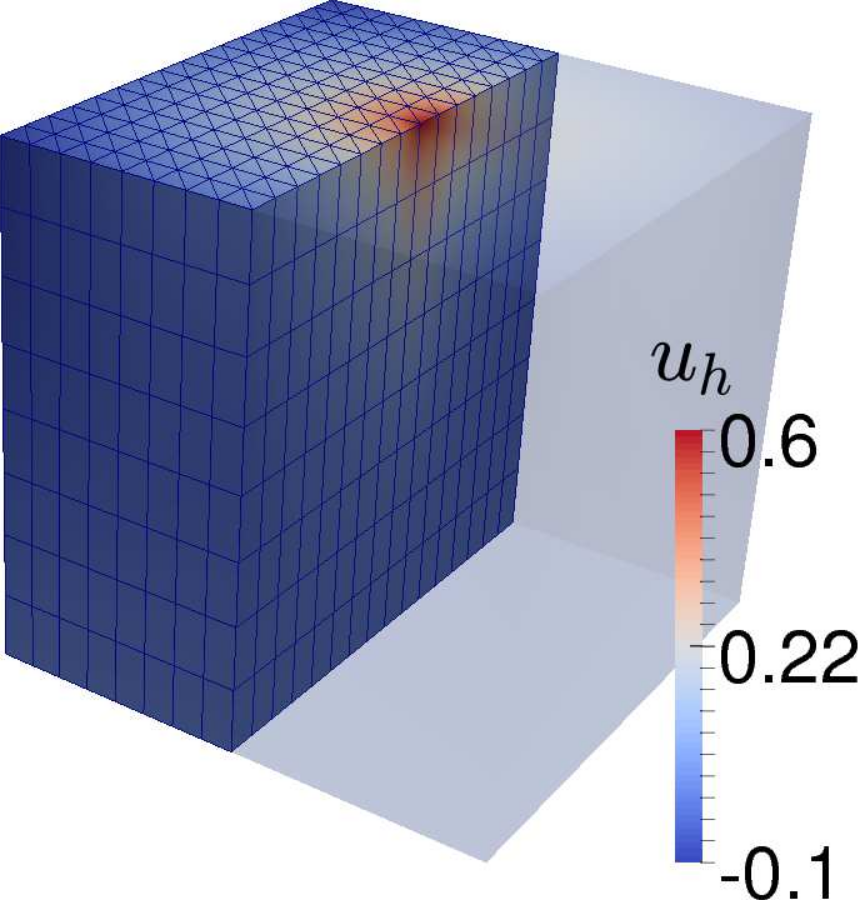}
        \caption{Total solution $u_h$}
    \end{subfigure}
    \caption{The splitting terms and solution itself found when solving \eqref{eq:model}-\eqref{eq:model-bc-d} on the unit domain with line source intensity $f(z)=z^3$ and the Dirichlet boundary data given in  \eqref{eq:exp1}.}
    \label{fig:infline-decomposition}
    
\vspace{3em}
    \captionsetup{labelformat=andtable}
    \caption{Convergence rates $p$ of the error in different norms, found when applying either the standard FE or SSB-FE method. The error is computed for the full domain $\Omega$ as well as for the domain with the singularity removed, i.e., $\Omega_R$ ($R=0.2$).}
    \addtocounter{figure}{-1}
    \begin{subfigure}{0.98\textwidth}
    \centering
     \caption{Convergence rates $p$ of the error $\norm{u-u_h}{}$ obtained when applying the standard FEM. }

\begin{tabular}{llrlrlrl}
\toprule
$h_{\Vert}$ & $h_\perp$ &  $L^2( \Omega)$ &  $p$ &  $L^2( \Omega_R)$ &  $p$ &  $H^1(\Omega_R)$ &  $p$ \\
\midrule
1/16 & 1/4 &         1.4e-2 &   &                    5.0e-3 &   &                    8.3e-2 &   \\
     & 1/8 &         6.0e-3 &  1.2 &                    1.4e-3 &  1.9 &                    4.4e-2 &  0.9 \\
     & 1/16 &         2.9e-3 &  1.1 &                    3.6e-4 &  1.9 &                    2.3e-2 &  1.0 \\
     & 1/32 &         1.5e-3 &  1.0 &                    1.5e-4 &  1.3 &                    1.4e-2 &  0.7 \\
     & 1/64 &         7.8e-4 &  0.9 &                    1.3e-4 &  0.2 &                    1.1e-2 &  0.3 \vspace{0.6em} \\ 
1/64 & 1/4 &         1.4e-2 &   &                    5.0e-3 &   &                    8.3e-2 &   \\
     & 1/8 &         6.0e-3 &  1.2 &                    1.4e-3 &  1.9 &                    4.3e-2 &  0.9 \\
     & 1/16 &         2.8e-3 &  1.1 &                    3.6e-4 &  1.9 &                    2.1e-2 &  1.0 \\
     & 1/32 &         1.4e-3 &  1.0 &                    1.0e-4 &  1.8 &                    1.0e-2 &  1.0 \\
     & 1/64 &         7.0e-4 &  1.0 &                    4.4e-5 &  1.2 &                    5.7e-3 &  0.9 \\
\bottomrule
\end{tabular}

\label{tab:infinite-exp1-u}
\vspace{2em}

    \caption{Convergence rates $p$ of the error $\norm{w-w_h}{}$ obtained when applying the SSB-FEM.}
\begin{tabular}{llrlrlrlrl}
\toprule
$h_{\Vert}$ & $h_\perp$ &  $L^2( \Omega)$ &  $p$ &  $H^1( \Omega)$ &  $p$ &  $L^2( \Omega_R)$ &  $p$ &  $H^1(\Omega_R)$ &  $p$ \\
\midrule
1/16 & 1/4 &         2.4e-2 &   &         2.4e-1 &   &                    1.9e-2 &   &                    2.2e-1 &   \\
     & 1/8 &         6.0e-3 &  2.0 &         1.2e-1 &  1.0 &                    4.5e-3 &  2.1 &                    1.0e-1 &  1.1 \\
     & 1/16 &         1.5e-3 &  2.0 &         5.9e-2 &  1.0 &                    1.1e-3 &  2.0 &                    4.9e-2 &  1.1 \\
     & 1/32 &         3.8e-4 &  2.0 &         2.9e-2 &  1.0 &                    2.8e-4 &  2.0 &                    2.4e-2 &  1.0 \\
     & 1/64 &         9.5e-5 &  2.0 &         1.5e-2 &  1.0 &                    7.1e-5 &  2.0 &                    1.2e-2 &  1.0 \vspace{0.6em} \\
1/64 & 1/4 &         2.4e-2 &   &         2.4e-1 &   &                    1.9e-2 &   &                    2.2e-1 &   \\
     & 1/8 &         6.0e-3 &  2.0 &         1.2e-1 &  1.0 &                    4.5e-3 &  2.1 &                    1.0e-1 &  1.1 \\
     & 1/16 &         1.5e-3 &  2.0 &         5.9e-2 &  1.0 &                    1.1e-3 &  2.0 &                    4.9e-2 &  1.1 \\
     & 1/32 &         3.8e-4 &  2.0 &         2.9e-2 &  1.0 &                    2.8e-4 &  2.0 &                    2.4e-2 &  1.0 \\
     & 1/64 &         9.4e-5 &  2.0 &         1.5e-2 &  1.0 &                    7.1e-5 &  2.0 &                    1.2e-2 &  1.0 \\
\bottomrule
\end{tabular}

\label{tab:infinite-exp1-w}
    \end{subfigure}
\end{figure*}

Let 
\begin{align*}
\Lambda= \lbrace (x_0,y_0,z) : z \in (0, 1) \rbrace
\end{align*}
be a line cutting vertically through the unit cube domain $\Omega$, with $x_0,y_0$ chosen as is illustrated in Figure \ref{fig:prism}. We solve \eqref{eq:model}-\eqref{eq:model-bc-d} with line intensity $f(z)=z^3$ and Dirichlet boundary conditions
\begin{align}
u_D = -\frac{1}{2\pi} \left( z^3 \ln(r) - 1.5 z \, r^2 \, (\ln({r})-1) \right) \text{ on } \partial \Omega,
\label{eq:exp1}
\end{align}
for which we have the analytic solution
\begin{align}
u_a = -\frac{1}{2\pi} \left( z^3 G  + w \right), \quad G = \ln(r), \quad w = - 1.5 z \, r^2 \, (\ln({r})-1).
\end{align}
The solution and its splitting terms are shown in Figure \ref{fig:infline-decomposition}.

Table \ref{tab:infinite-exp1-u} shows the numerical results obtained when solving for $u_h$ directly by \eqref{eq:varform-uh}. Examining the $L^2(\Omega)$-error of the solution, we see that the solution converges in this norm sub-optimally, i.e., with only order $h^{1-\epsilon}$ for some small $\epsilon>0$. The convergence order improves to order $h^2$, i.e., optimal, when a small region is removed around the line. This agrees with the results in \cite{koppl2014}, where quasi-optimal convergence was proven for the point source in a 2D domain as long as a small area was removed around the point source. The $H^1$-error for the entire domain is not given as the error does not converge in this norm. 

The results in Table \ref{tab:infinite-exp1-u} are also consistent with our observation that $u$ exhibits anisotropic regularity. Comparing the errors for $h_\Vert=1/16$ and $h_\Vert=1/64$, one can conclude that a refinement along with the line makes close to no difference in the $L^2$-error taken on the full domain. It follows that the perpendicular error dominates the parallel error; this is consistent with our observation in Section \ref{sec:decomposition} that the solution exhibits a regularity loss in the plane perpendicular to the line. In contrast, both $L^2$ and $H^1$-errors are affected by the vertical refinement when the cylinder $U_{R=0.2}$ is removed from the domain. This is to be expected as we know from Section \ref{sec:decomposition} that the solution exhibits high and isotropic regularity as long as a small volume is removed around the line. 

Table \ref{tab:infinite-exp1-w} shows the numerical results obtained when solving for $w_h$ by \eqref{eq:varform-wh} and comparing with the analytical expression for $w$. As the correction term here solves the Poisson equation with right-hand side $f''(z) \ln(r) = z \ln(r) \in L^2(\Omega)$, we see here optimal convergence rates. The error is not affected by removing a small region around the line; this makes sense as $w$ is $H^2$ regular in the entire domain. We also here see evidence of solution anisotropy, as the error is comparable for $h_\Vert=1/16$ and $h_\Vert=1/64$. This is natural as $w$ adopts the regularity of $f$ in the vertical direction, meaning that the correction function $w$ for this test case is smooth with respect to $z$. The table also lists the convergence properties of the error in the $H^1$-norm, which was similarly found to be of optimal order.

\subsection{Convergence test for arbitrary line source}
\label{sec:num-finite-line}

\begin{figure*}
    \centering
    \begin{subfigure}[t]{0.32\textwidth}
        \centering
        \includegraphics[height=1.5in]{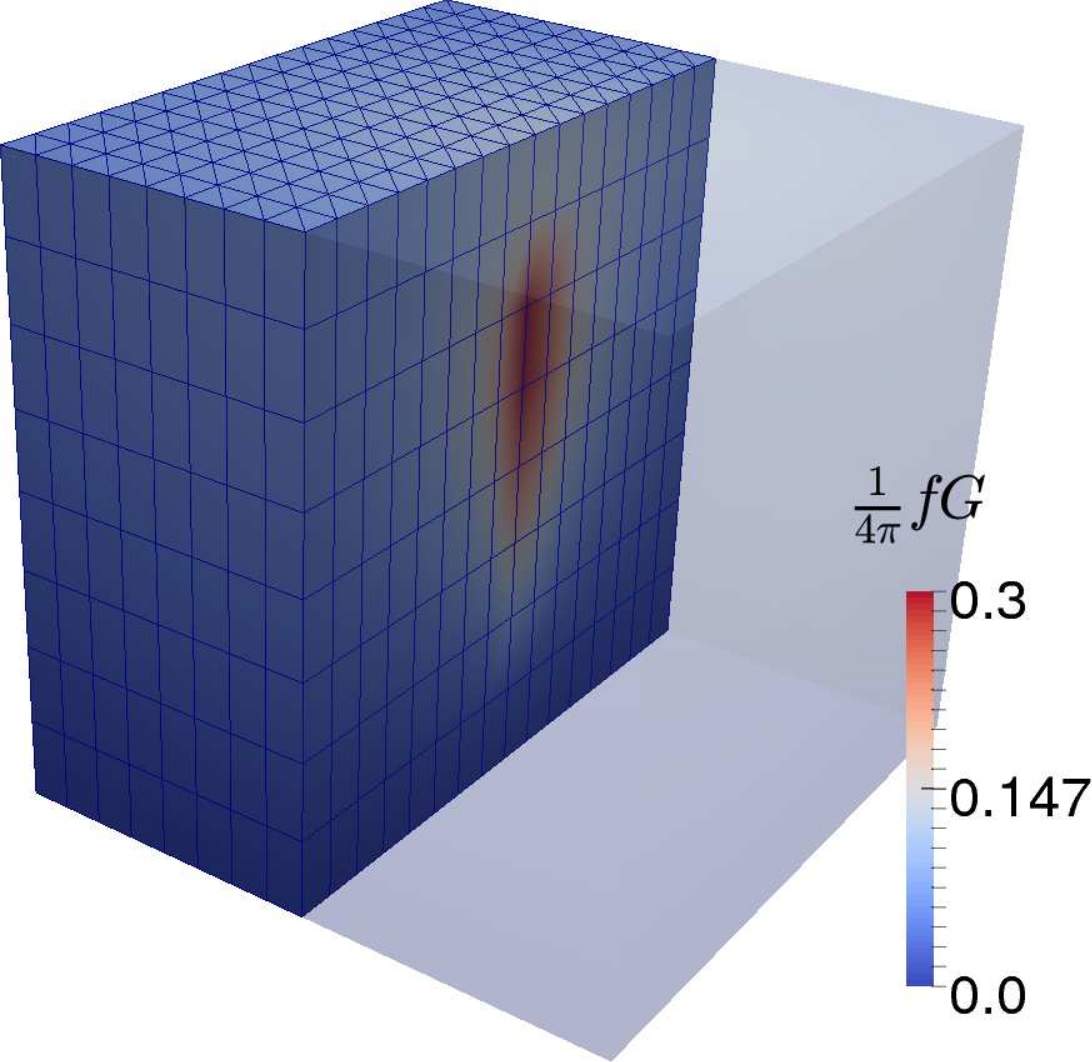}
        \caption{$ \frac{1}{4\pi} f(z) G $}
    \end{subfigure}
    \begin{subfigure}[t]{0.32\textwidth}
        \centering
        \includegraphics[height=1.5in]{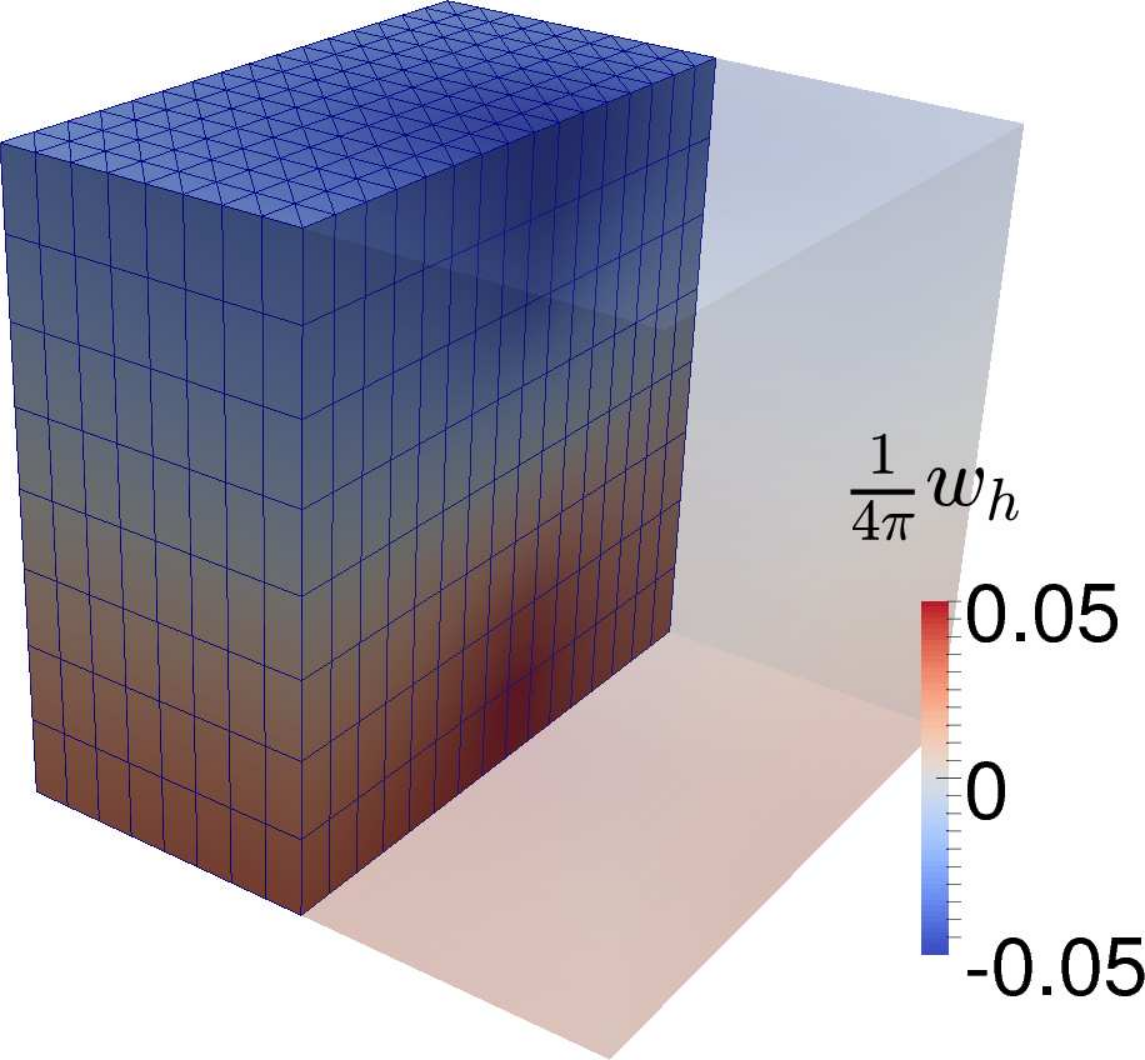}
        \caption{ Correction $ \frac{1}{4\pi} w_h $}
    \end{subfigure}
    \begin{subfigure}[t]{0.32\textwidth}
        \centering
        \includegraphics[height=1.5in]{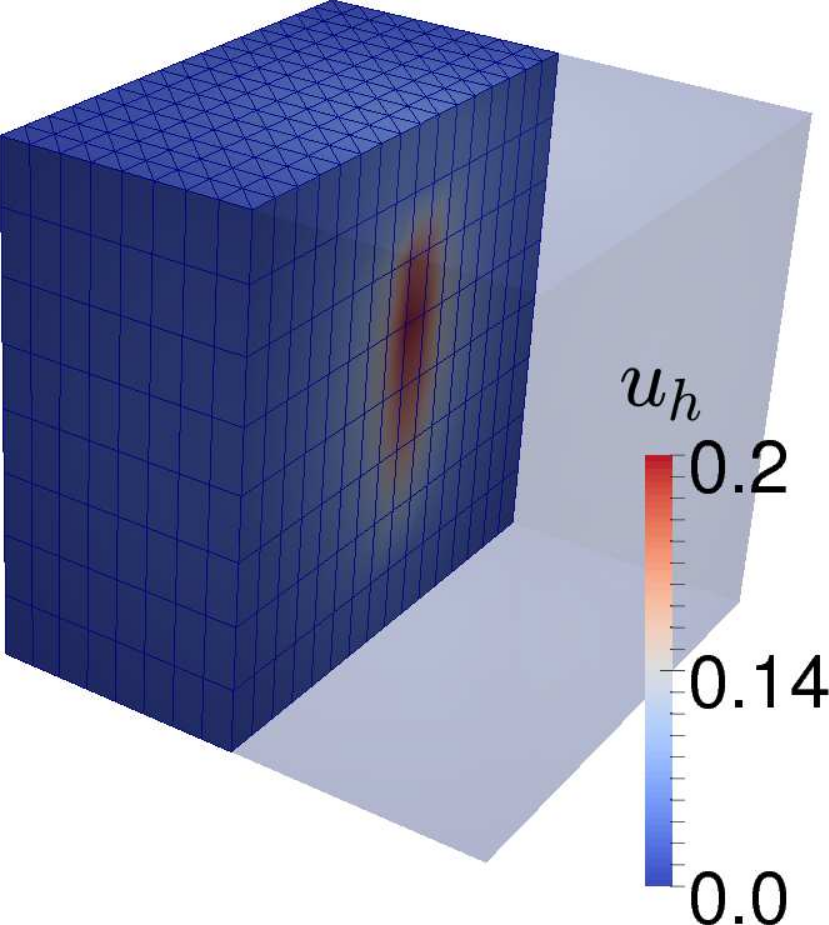}
        \caption{Total solution $u_h$}
    \end{subfigure}
      \caption{The splitting terms and solution itself found when solving \eqref{eq:model}-\eqref{eq:model-bc-d} on the unit cube domain with line source intensity $f(z)=z$ and the Dirichlet boundary data in \eqref{eq:exp1}.}
    \label{fig:finite-line-decomposition}
    
    \vspace{2em}
    
    \captionsetup{labelformat=andtable}
   \caption{Convergence rates $p$ of the error in different norms found when applying either the standard FE method or the SSB-FE method. The error is computed for the full domain $\Omega$ as well as for the domain with the singularity removed, i.e., $\Omega_R$ ($R=0.2$).}
    \addtocounter{figure}{-1}
    \begin{subfigure}{1.0\textwidth}
    \centering
    \caption{Convergence rates $p$ of the error $\norm{u-u_h}{}$ obtained when applying the standard FEM.}
\begin{tabular}{llrrrrrrrr}
\toprule
$h_{\Vert}$ & $h_\perp$ &  $L^2( \Omega)$ &  $p$ &  $L^2( \Omega_R)$ &  $p$ &  $H^1(\Omega_R)$ &  $p$ \\
\midrule
1/16  & 1/4  &         1.0e-2 &        &                          6.5e-3 &      &                    1.4e-1 &      \\
      & 1/8  &         5.8e-3 & 0.8 &                          2.2e-3 &  1.6 &                    6.3e-2 & 1.2 \\
      & 1/16 &         3.2e-3 & 0.9 &                          6.4e-4 &  1.7 &                    3.3e-2 & 0.9 \\
      & 1/32 &         1.8e-3 & 0.9 &                          4.5e-4 &  0.5 &                    1.6e-2 & 1.0 \\
      & 1/64 &         1.1e-3 & 0.7 &                          4.6e-4 &      &                    9.8e-3 & 0.7 \\ \\
1/128 & 1/4  &         1.0e-2 &        &                          6.5e-3 &      &                    1.4e-1 &      \\
      & 1/8  &         5.6e-3 & 0.9 &                          2.2e-3 &  1.6 &                    6.4e-2 & 1.1 \\
      & 1/16 &         3.0e-3 & 0.9 &                          5.5e-4 &  2.0 &                    3.3e-2 & 1.0 \\
      & 1/32 &         1.6e-3 & 0.9 &                          1.8e-4 &  1.6 &                    1.5e-2 & 1.1 \\
      & 1/64 &         8.3e-4 & 0.9 &                          1.3e-4 &  0.4 &                    8.0e-3 & 0.9 \\
\bottomrule
\end{tabular}

\label{tab:finite-exp2-u}

\vspace{2em}

    \caption{Convergence rates $p$ of the error $\norm{w-w_h}{}$ obtained when applying the SSB-FEM.}\begin{tabular}{llrrrrrrrr}
\toprule
$h_{\Vert}$ & $h_\perp$ &  $L^2( \Omega)$ &  $p$ &  $H^1( \Omega)$ &  $p$ &  $L^2( \Omega_R)$ &  $p$ &  $H^1( \Omega_R)$ &  $p$ \\
\midrule
1/16  & 1/4  &         9.7e-4 &     &         1.4e-2 &        &                    7.5e-4 &     &                    1.2e-2 &      \\
      & 1/8  &         2.9e-4 & 1.7 &         7.8e-3 & 0.8    &                    2.0e-4 & 1.9 &                    6.2e-3 & 1.0 \\
      & 1/16 &         1.1e-4 & 1.4 &         4.9e-3 & 0.7   &                    7.9e-5 & 1.4 &                    4.0e-3 & 0.6 \\
      & 1/32 &         7.4e-5 & 0.6 &         3.8e-3 & 0.4 &                    6.0e-5 & 0.4 &                    3.2e-3 & 0.3 \\
      & 1/64 &         6.7e-5 & 0.1 &         3.4e-3 & 0.2 &                    5.7e-5 & 0.1 &                    3.0e-3 & 0.0 \\ \\
1/128 & 1/4  &         9.5e-4 &     &         1.3e-2 &        &                    7.4e-4 &     &                    1.2e-2 &      \\
      & 1/8  &         2.7e-4 & 1.8 &         7.1e-3 & 0.9 &                    1.8e-4 & 2.0 &                    5.5e-3 & 1.1 \\
      & 1/16 &         7.3e-5 & 1.9 &         3.7e-3 & 0.9 &                    4.3e-5 & 2.1 &                    2.8e-3 & 1.0 \\
      & 1/32 &         1.9e-5 & 1.9 &         1.9e-3 & 0.9 &                    1.1e-5 & 2.0 &                    1.4e-3 & 1.0 \\
      & 1/64 &         5.1e-6 & 1.9 &         1.0e-3 & 0.9 &                    3.0e-6 & 1.9 &                    7.9e-4 & 0.9 \\
\bottomrule
\end{tabular}

\label{tab:finite-exp2-w}
    \end{subfigure}
\end{figure*}

We again solve \eqref{eq:model}-\eqref{eq:model-bc-d} on the unit cube, with $\Lambda$ as in Section \ref{sec:num-infinite-line}, with the exception that it now has two endpoints contained in the interior of $\Omega$:
\begin{align}
\Lambda = \lbrace (x,y,z) \in \Omega : x=x_0, y=y_0, z \in (0.2,0.8) \rbrace.
\end{align}
The line source is prescribed an intensity $f(z) = z$ and Dirichlet boundary conditions are assigned on the boundary, 
\begin{align}
u_D &= \frac{1}{4\pi} \left(  z \ln( \frac{\sqrt{r^2+(z-b)^2} - (z-b)}{\sqrt{r^2+(z-a)^2} - (z-a)} \right) \\
&  - \left( \sqrt{(z-a)^2+r^2} - \sqrt{(z-b)^2+r^2}  \right) \text{ on } \partial \Omega. 
\end{align}
The solution to \eqref{eq:model}-\eqref{eq:model-bc-d} then admits an analytic solution, namely
\begin{align}
u = \frac{1}{4\pi} \left( E(f) G + w \right), \quad E(f) = z,  \quad w = &  -\sqrt{(z-a)^2+r^2} + \sqrt{(z-b)^2+r^2},
\end{align}
where $G$ is given by \eqref{eq:G}. The solution, along with its splitting terms, is shown in Figure \ref{fig:finite-line-decomposition}.

Table \ref{tab:finite-exp2-u} shows the numerical results obtained when solving for $u_h$ by the standard FE method. Examining again the $L^2(\Omega)$-error of the solution, we see that the solution converges in this norm sub-optimally, i.e., with only order $h^{1-\epsilon}$ for some small $\epsilon>0$. The convergence order becomes optimal, i.e., of order $h^2$, when a small region is removed around the line. 

As in the previous example, the $L^2$-error shows the anisotropic nature of the solution, with an increase in the vertical refinement from $h_\Vert=1/16$ to $h_\Vert=1/128$ leading only to a small reduction in the error. It is thus clear that the error for this test case is also dominated by the horizontal error. Removing a small region around the line again removes the anisotropic nature of the error; this is evident from the fact that the convergence rates of the error in $L^2(\Omega_{R})$ and $H^1(\Omega_{R})$-norms quickly deteriorate when the vertical refinement becomes coarser than the horizontal refinement. 

Table \ref{tab:finite-exp2-w} shows the numerical results obtained when solving for $w_h$ with the SSB-FE method and comparing the solution with the analytical expression for $w$. The results are again better when compared to solving for $u_h$ directly, with optimal order convergence observed in both $L^2$ and $H^1$-norms. The convergence rates remain nearly unaltered when a small volume is removed around, confirming our observation that $w$ exhibits high-regularity on the entire domain $\Omega$.

\subsection{Vascular network of the brain}
\label{sec:num-brain}
Finally, let us demonstrate the capability of the SSB-FE method in handling problems with a large number of line segments. To this end, we solve the line source problem with both Dirichlet and Neumann type boundary conditions, 
\begin{subequations}
\begin{align}
 - \Delta u &= f \delta_\Lambda  \quad   &\text{in} \, \Omega, \label{eq:model-full}\\
 u &= u_D  \quad   &\text{on} \, \partial \Omega_D, \label{eq:model-full-bc-d}\\ 
 \nabla u\cdot \mathbf{n} &= 0  \quad   &\text{on} \, \partial \Omega_N, \label{eq:model-full-bc-n}
\end{align}
\end{subequations}
where $\mathbf{n}$ denotes the outward directed boundary normal of the domain, $u_D$ some given boundary data, and $\partial \Omega_D$ and $\partial \Omega_N$ the Dirichlet and Neumann boundaries respectively, assumed to satisfy $\partial \Omega_D \cap \partial \Omega_N = \emptyset$ and $\partial \Omega_D \cup \partial \Omega_N = \partial \Omega$. The line sources are taken from a graph describing the vascular network in a human brain \cite{brain}, the same dataset as was illustrated in Figure \ref{fig:drawing}. This dataset consists of $n_a=264$ and $n_v = 2606$ line segments describing the arterial and venous systems of the brain, respectively. The arterial system is responsible for providing the brain with oxygenated blood; its counterpart is the venous system, which is tasked with returning de-oxygenated blood back to the heart.

The process of oxygen delivery to biological tissue by micro-circulation is commonly modelled by means of coupled 1D-3D flow models, where the graph is seen as a 1D network and endowed with its own flow equation. The 1D flow equation can then be coupled to the 3D flow equation for the tissue by Starling's law of filtration. In this work, however, we consider only the 3D flow equation with a collection of line sources. For this reason, it is necessary to prescribe a certain mass flux to each line segment. For the test case presented here, the linear mass flux from each vessel $i$ was taken proportional to its (given) vessel radius $R_i$:
\begin{align}
f_{i}(t) = \gamma_i R_{i} \quad \text{for } t\in (0, L_i),
\label{eq:brain-f}
\end{align}
where the proportionality constant is taken as $\gamma_i=1.0$ for arterial blood vessels and $\gamma_i=-0.1$ for venous blood vessels. The graph was then embedded in a spherical domain loosely representing the skull. A Dirichlet pressure was assigned to the bottom portion of the domain,
\begin{align}
u_D = 1 \text{ for } \mathbf{x} \in \partial \Omega_D, \quad
\partial \Omega_D = \lbrace (x,y,z) \in \partial \Omega : x<H \rbrace,
\label{eq:brain-bc-D}
\end{align}
and no-flow out of the remaining boundary,
\begin{align}
\nabla u \cdot \mathbf{n} = 0  \text{ for } \mathbf{x} \in \partial \Omega_N, \quad \partial \Omega_N = \lbrace (x,y,z) \in \partial \Omega : x > H \rbrace.
\label{eq:brain-bc-N}
\end{align}
The solution then admits a splitting 
\begin{align}
u = \frac{1}{4\pi} \sum_{i=1}^{n_a+n_v} \left( E_i(f) G + w  \right)
\end{align}
where the extension is given by 
\begin{align}
E_i(f)(\mathbf{x}) = \gamma_i R_{i} \text{ for } \mathbf{x} \in \Omega,
\end{align}
and $w$ solves the boundary value problem Laplace equation (with right-hand side $F=0$) with boundary conditions 
\begin{align}
w_D &= 4 \pi u_D - \sum_{i=1}^{n_a + n_v} E_{i}(f) G_{i} &\text{ for } \mathbf{x} \in \partial \Omega_D, \\
\nabla w \cdot \mathbf{n} &= - \sum_{i=1}^{n_a + n_v} \nabla E_{i}(f) G_{i} \cdot \mathbf{n}  &\text{ for } \mathbf{x} \in \partial \Omega_N.
\end{align}
We see that for the choice of piecewise constant linear mass flux from each line segment, the singularity subtraction technique reduces the line source problem to a boundary value problem. 

\begin{figure*}
    \centering
    \begin{subfigure}[t]{0.99\textwidth}
        \centering
         \includegraphics[height=3.3in]{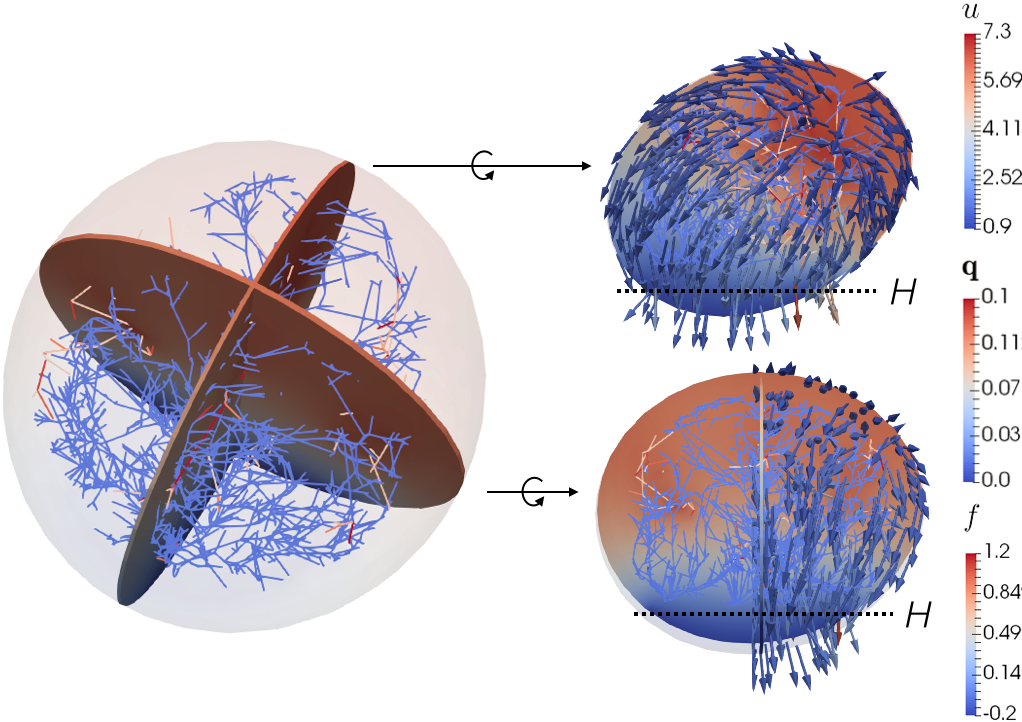}
  \caption{(Left) perspective view of the brain, with two cross section views showing the pressure in the plane dividing the left/right hemisphere and front/rear hemisphere, respectively. (Right top) side view and (right bottom) rear view of these two cross sections. No-flow was imposed in the region above $H$, while a a Dirichlet condition for the pressure was imposed in the region below $H$.}
  \label{fig:brain-bcs}
    \end{subfigure}
    
    \begin{subfigure}[t]{0.99\textwidth}
        \centering
           \includegraphics[height=3.0in]{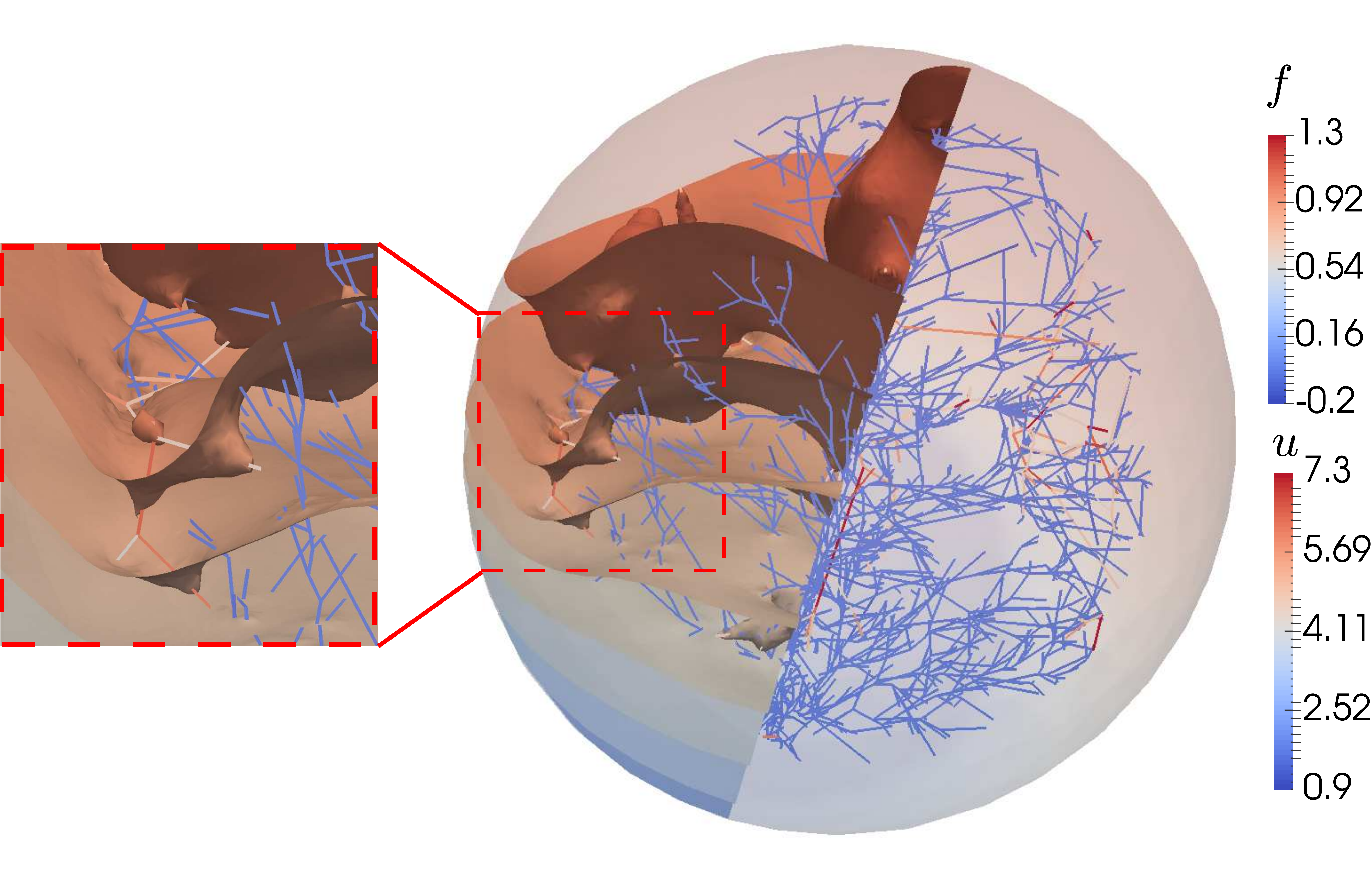}
           \caption{Contour plot of the pressure. The singular parts of the solution were interpolated using $k=4$; this was found necessary in order to resolve the logarithmic nature of the solution around the individual blood vessels.}
           \label{fig:brain-contour}
    \end{subfigure}
        \caption{Pressure and flux found when using the SSB-FEM to solve \eqref{eq:model-full}-\eqref{eq:model-full-bc-n} on a dataset describing the vascular system of the brain, with boundary conditions and source intensities set as specified in \eqref{eq:brain-f}-\eqref{eq:brain-bc-N}. The correction function $w$ was solved using linear elements ($k=1$) and the singular parts of the solution interpolated using quartic elements ($k=4$).}
    \label{fig:brain}
    
\end{figure*}

The results from using the SSB-FE method to solve this model problem are shown in Figure \ref{fig:brain}. Figure \ref{fig:brain-bcs} shows the source intensity $f_i$ prescribed for each edge in the graph and the resulting pressure solution $u$. The flux, found by a projection of $-\nabla u$, is also shown, from which we see that the solution satisfies the no-flow boundary condition imposed on $\partial \Omega_N$. As was discussed in Section \ref{sec:decomposition}, the line sources in the right-hand side of \eqref{eq:model-full}-\eqref{eq:model-full-bc-n} induce a logarithmic type singularity in the pressure profile around each line segment. Figure \ref{fig:brain-contour} shows a contour plot of the pressure, with a close-up view given of one region. The singular parts of the solution, i.e., the $E_i(f) G_i$-terms, were here discretized by interpolation onto a fine mesh. In this close-up view the logarithmic profile of the solution is indeed visible, showing that the SSB-FE method is capable of producing good discrete approximations of the pressure. The logarithmic behaviour of the solution is especially pronounced around line segments representing arteries; this is to be expected as the arteries were prescribed a comparatively high mass flux. 

\begin{figure}
\centering
       \includegraphics[height=2.3in]{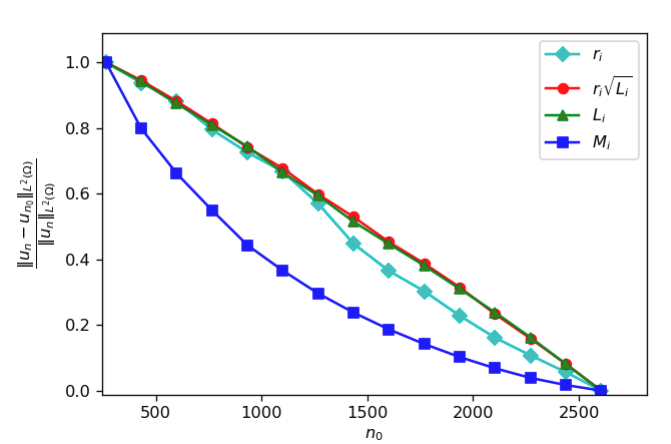}
      \caption{The model reduction error found when the graphs $\Lambda_{a}$ and $\Lambda_{v}$ were sorted montonically decreasing by vessel radius $r_i$, vessel length $L_i$, and the model reduction error $M_i$.}
       \label{fig:model-reduction}
\end{figure}

Lastly, let us consider the robustness of the solution with
respect to removing branches of the network. This serves as a proxy for more detailed studies regarding imperfect data acquisition and image segmentation, as well as for assessing changes in the actual arterial or venous network associated with diseases such as stroke. In particular, we ask the question as to how sensitive the solution (here represented by the $L^2$ norm of the error) is to a loss of a certain fraction of the edges of the graph.

Recall from Section \ref{sec:numerical-disc} that the modelling error associated with neglecting line segments $i=n_0,...,n_a+n_v$ is given by
\begin{align}
\norm{u-u_{n_0}}{L^2(\Omega)} 
& \lesssim \sum_{i={n_0}}^{n_a+n_v} \norm{E_i(f) G_i }{L^2(\Omega)} +  \underbrace{\norm{F_i}{L^2(\Omega)}}_{=0}  +\norm{w_D}{L^2(\Omega)} \\
& \lesssim  \sum_{i={n_0}}^{n_a+n_v} \norm{E_i(f) G_i }{L^2(\Omega)} + \norm{E_i(f) G_i }{H^\frac{1}{2}(\partial \Omega)} \\
&= \sum_{i={n_0}}^{{n_a+n_v}} \gamma_i R_i \left( \norm{G_i }{L^2(\Omega)} + \norm{G_i }{H^\frac{1}{2}(\Omega)} \right),
\label{eq:model-red-error}
\end{align}
where $F_i=0$ as the extensions $E_i(f)$ were chosen as constant functions.  This suggests that the modelling error associated with neglecting line segment $i$ depends the vessel radius $R_i$ as well as logarithmic term $G_i$, where the exact interplay between the two is not known. For this reason, we consider now four different enumerations of the edge segments. The two first enumerations were based on sorting monotonically decreasing by edge length $L_i$ and parametrized radius $r_i$, respectively. The third is formed by sorting monotonically decreasing with respect to $R_i \sqrt{L}_i$.
This was motivated by considering the $L^2$ norm of $G$ in the cylindrical domain $U_{R,H} = \lbrace  (r, \theta, z) : 0<r<R, \, 0<\theta<2\pi, \, -H < z < H \rbrace$ when $G$ corresponds to the line segment $\Lambda_i = \lbrace (r, \theta, z) \in \Omega_{R, H} : -L<z<L \rbrace$. We then have the approximation
\begin{align}
G &= \ln( \frac{\sqrt{r^2+(z+L)^2}-(z+L)}{\sqrt{r^2+(z-L)^2}-(z-L)} ) = \ln( \frac{\sqrt{1+(\frac{r}{z+L})^2}-1}{\sqrt{1+(\frac{r}{z-L})^2}-1} ) \\
& \approx \ln( \frac{r^2}{2(z-L)^2} ) - \ln(\frac{r^2}{2(z+L)^2}) 
= \frac{1}{2}  \left( \ln(z+L) - \ln(z-L) \right), \label{eq:approx}
\end{align}
where we have used the binomial expansion $\sqrt{1+x}=1+1/2x$. From this approximation, the $L^2$-norm of $G$ can be calculated explicitly, and we find that it scales $\norm{G}{L^2(U_{R,H})} \sim \sqrt{L}$. 

The fourth and final enumeration was made by numerically evaluating the $L^2$-norm of $G_i$ in the domain $\Omega$, and sorting monotonically decreasing by the approximation $M_i$ of the full modelling error:
\begin{align}
M_i &= {\gamma}_{i} R_i \norm{ G_{i} }{L^2(\Omega)}.
\label{eq:model-error-brain}
\end{align}

Figure \ref{fig:model-reduction} shows the actual (computed) error associated by removing a fraction of the edges based on the four enumerations given above. We see that there is essentially no correlation between the model error and segment length or radius  (as reflected in nearly straight lines in the figure). This suggests that the modelling error in \eqref{eq:model-red-error} is difficult to quantify in terms of the model parameters, even for the simple test case presented here. The $L^2$-norm of each $G_i$-term depends on the line segment length as well as its location within the domain. For this reason, the approximation shown in \eqref{eq:approx} is not helpful as $\norm{G_i}{L^2(\Omega)}$ can depend strongly on the location of $\Lambda_i$ relative to the domain. The enumeration using the model reduction error estimate, $M_i$, is the most successful of the four; Figure \ref{fig:model-reduction} shows that the reduction error decreases near quadratically to zero when the dataset is sorted using this parameter. This suggest that, to assess the impact of removing a line segment $\Lambda_i$, it is appropriate to calculate the norm $\norm{G_i}{L^2(\Omega)}$ numerically. 

Although this study is beholden to several academic simplifications, observations of this nature have several practical implications. Firstly, imaging technology based on finite resolution will primarily omit edges with small radius. Figure \ref{fig:model-reduction} implies that it is not appropriate to assume that these edges can be neglected, as a smaller edge radius was not seen to imply negligible impact. For applications, it may therefore be desirable to consider data sets which are augmented from the registered data with synthetic vessels of finer radius. Secondly, these results also give an understanding of the associated risks of failure of blood flow to the brain. Indeed, if the $L^2$-error of pressure is taken as a proxy for the change in oxygen perfusion in the brain, we again see that the term $M_i$ provides an indicator for the most sensitive edges of the arterial and venous systems. While these observations are encouraging, we emphasize that a more detailed study is required in order to claim medical relevance. Such a study will need to include (among other things) a realistic geometry for not just of the arterial and venous network, but also the brain tissue; a coupling between the source terms of the edges with transport on the edge network (coupled 1D-3D); a time dependent right-hand side to account for the pulsating nature of the macro-vascular system, as the assumption of stationary flow is generally only valid for micro-circulation (see e.g. the model in \cite{secomb2010}); as well as an assessment of the non-linearities associated with non-Newtonian fluid flow.

\textbf{Acknowledgements}
The authors thank J. Reichenbach and A. Deistung for bringing our attention to the data used in section \ref{sec:num-brain} \cite{brain}, and E. Hanson and E. Hodneland for providing us with the data segmentation and tree extraction.

\section{Conclusions}
We studied an elliptic equation having line sources in a 3D domain. The line sources act as Dirac measure defined on a line causing the solutions to be singular on the line itself. Central to this work is the result that the solution admits a split into a singular and a regular part. This allows us to study the nature of the solution as well as to develop a numerical algorithm for solving the problem.  Mathematically, we see that the solution has anisotropic regularity, it is smooth along the line source and the line singularity acts as a Dirac point measure in a 2D domain. Our numerical approach solves for the regular part only and therefore obtains optimal convergence rates.  We illustrate our approach for several numerical examples including a data set describing the vascular system of a human brain. Our solution approach is mesh-independent and can be adapted to a variety of discretizations. 
\label{sec:conclusions}

\bibliographystyle{apalike}

\end{document}